\input epsf.tex % To insert eps figures.

% Some font names:
\font\titlefont=cmr12 at 14pt
\font\sectionfont=cmbx10 at 12 pt
\font\curly=eusm10 at 10.7pt
\font\Bbb=msbm10
\font\smBbb=msbm7
%\font\subsectfont=cmss10 at 11pt
\font\subsectfont=cmbx10 at 11pt

\def\subsect#1\par{\bigskip \leftline{\subsectfont #1 } \nobreak\smallskip}

% Some typesetting parameters:
\magnification=\magstep1
\baselineskip=15pt
\mathsurround=1pt
\abovedisplayskip=8pt plus 3pt minus5pt
\belowdisplayskip=8pt plus 3pt minus5pt
\fontdimen16\textfont2=2.5pt  %subscript, no superscript
\fontdimen17\textfont2=2.5pt  %subscript, with superscript
\fontdimen14\textfont2=4.5pt  %superscript, textstyle
\fontdimen13\textfont2=4.5pt  %superscript, displaystyle

\def\page{\vfill\break}
\def\pf #1. {\noindent {\it #1}.\enskip}

\def\Diff{\hbox{\rm Diff}}
\def\Emb{\hbox{\rm Emb}}
\def\Aut{\hbox{\rm Aut}}

\def\int{\hbox{\rm int}}
\def\qed{\quad\hfill \rlap{$\sqcup$}$\sqcap$\par\medskip\smallskip} 

\def\:{\,\colon} % Used for function notation.

\def\scrptA{\lower.2pt\hbox{\curly A}}
\def\scrptC{\lower.2pt\hbox{\curly C}}
\def\scrptE{\lower.2pt\hbox{\curly E}}
\def\scrptG{\lower.2pt\hbox{\curly G}}
\def\scrptH{\lower.2pt\hbox{\curly H}}
\def\scrptK{\lower.2pt\hbox{\curly K}}
\def\scrptL{\lower.2pt\hbox{\curly L}}
\def\scrptM{\lower.2pt\hbox{\curly M}}
\def\scrptT{\lower.2pt\hbox{\curly T}}

\def\C{\scrptC}
\def\E{\scrptE}

\def\M{\scrptM}

\def\Q{\hbox{\Bbb Q}}
\def\R{\hbox{\Bbb R}}
\def\Z{\hbox{\Bbb Z}} 
\def\smQ{\hbox{\smBbb Q}}

\def\CP{\hbox{\Bbb C}{\rm P}}

\def\Cn{\C^n}
\def\Cnn{\C^{n,n}}
\def\Cnk{\C^{n,k}}
\def\Cinfone{\C^{\infty,1}}
\def\Cnkplusone{\C^{n,k+1}}
\def\Cnzero{\C^{n,0}}
\def\Cnone{\C^{n,1}}
\def\Cntwo{\C^{n,2}}
\def\Mnk{\M^{n,k}}
\def\Mn{\M^n}
\def\Mnzero{\M^{n,0}}

\def\Minf{\M^\infty}

\def\phi{\varphi}

\def\bdy{\partial}
\def\tilde{\widetilde}

\def\incl{\hookrightarrow}
\def\To{\longrightarrow}

\def\mapright#1{\smash{\mathop{\ \longrightarrow\ }\limits^{#1}}}
\def\mapleft#1{\smash{\mathop{\ \longleftarrow\ }\limits^{#1}}}
\def\iso{\cong}

\outer\def\section#1\par{\vskip0pt plus.05\vsize\penalty-250\vskip0pt plus-.05\vsize 
\vskip1.2cm
\leftline{\sectionfont#1}\medskip}

%\def\subsection #1\par{\noindent {\bf#1}\smallskip}

%%%%%%%%%%%%%%%%%%%%%%%%%%%%%%%%%%%%%%%%

\vbox{}\vskip 4mm
\centerline{\titlefont A Short Exposition of the Madsen-Weiss Theorem}

\vskip3mm
\centerline{Allen Hatcher}

\vskip18mm

The theorem of Madsen and Weiss [MW] identifies the homology of mapping class groups of surfaces, in a stable dimension range, with the homology of a certain infinite loopspace.  This result is not only intrinsically interesting, showing that two objects that appear to be quite different turn out to be homologically equivalent, but it also allows explicit calculations of the stable homology, rather easily for rational coefficients as we show in Appendix C, and with considerably more work for mod $p$ coefficients as done by Galatius in [G2].  Outside the stable dimension range the homology of mapping class groups appears to be quite complicated and is still very poorly understood, so it is surprising that there is such a simple and appealing description in the stable range.  

The Madsen-Weiss theorem has a very classical flavor, and in retrospect it seems that it could have been proved in the 1970s or 1980s since the main ingredients were available then.  However, at that time it was regarded as very unlikely that the stable homology of mapping class groups could be that of an infinite loopspace.  The initial breakthrough came in a 1997 paper of Tillmann [T] where this unexpected result was proved.  It remained then to determine whether the infinite loopspace was a familiar one.  A conjecture in this direction was made in a 2001 paper of Madsen and Tillmann [MT], with some supporting evidence, and this conjecture became the Madsen-Weiss theorem.

The original proof of the Madsen-Weiss theorem was rather lengthy, but major simplifications have been found since then.  The purpose of the present paper is to present a proof that uses a number of these later simplifications, particularly some due to Galatius and Randal-Williams which make the proof really quite elementary, apart from the three main classical ingredients:

\smallskip
\item{(1)} The fact that the diffeomorphism group of a closed orientable surface has contractible components once the genus of the surface is at least $2$.  This is originally a theorem of Earle-Eells [EE] proved by analytic methods, but a purely topological proof was found soon after by Gramain [Gr].  Also needed is the extension of this result to compact orientable surfaces with boundary, originally a result in [ES], but proved in the Gramain paper as well. 
In fact Gramain's proof in both cases is quite simple so we present it in Appendix~B.

\item{(2)} Harer stability [Har], the fact that the $i$th homology group of the mapping class group of a compact orientable surface is independent of the genus and the number of boundary components, once the genus is sufficiently large with respect to $i$.  Subsequent improvements in the stable dimension range were made by Ivanov [I], Boldsen [B], and Randal-Williams [RW], and a significant gap in Harer's proof was filled by Wahl [W1].  For a full exposition of the current state of the art on this theorem see [W2]. As yet no really simple proof has been found.

\item{(3)} The Group Completion Theorem from around 1970.  A nice exposition of this fundamental result in algebraic topology was published in [MS], and an illuminating general discussion can be found in [A].  The original version seems to be due to Barratt and Priddy [BP].  This theorem seems not to have found its way into textbooks yet, so we give a textbook-style exposition of it in Appendix~D following an argument told to us by Galatius.  The only ingredients for this are results already available in textbooks.

\smallskip\noindent
The strategy we follow for proving the Madsen-Weiss theorem, just as in [GRW], is to use (3) to prove a theorem that is independent of (1) and (2).  This theorem can be phrased as saying that the classifying space of the group of compactly supported diffeomorphisms of an infinite-genus surface is homologically equivalent to a certain infinite loopspace (or more properly, one component of this infinite loopspace).  Then one can quote (1) and (2) to restate this result in terms of the homology of mapping class groups of compact surfaces in the stable range.

A large part of the proof works for manifolds of arbitrary dimension, not just surfaces, so we present it in that generality.  In the case of zero-dimensional manifolds the whole proof goes through and in fact is considerably simpler than for surfaces, yielding the Barratt-Priddy-Quillen theorem that the infinite symmetric group $\Sigma_\infty = \cup_n\Sigma_n$ has the same homology as one component of $\Omega^\infty S^\infty$.  

Here is a quick outline of the paper.  In section 1 we present basic definitions and the motivating constructions leading to a precise statement of the Madsen-Weiss theorem.  The proof of the theorem is contained in section~2, with the two hardest steps postponed to sections~3 and~4.  In section~5 we give extensions to the cases of nonorientable surfaces and surfaces with punctures, along with the Barratt-Priddy-Quillen theorem and its analog for braid groups due to F.\ Cohen.  After this we have four appendices:  Appendix~A giving basic information on classifying spaces for diffeomorphism groups, Appendix~B proving the Earle-Eells theorem following the method of Cerf and Gramain, Appendix~C giving the calculation of the stable rational homology in both the orientable and nonorientable cases, and finally Appendix~D proving the Group Completion Theorem.

\medskip\noindent
{\bf Acknowledgments:} \enskip Thanks are due to Tom Church, S\o ren Galatius, Oscar Randal-Williams, and Ulrike Tillmann for their comments on an early version of this exposition.

\section 1. The Scanning Map

The statement of the Madsen-Weiss theorem involves a certain infinite loopspace, so we will begin by describing what this space is and how it arises.  The root of the idea can be traced back a long way, to the Pontryagin-Thom construction around 1950.  Another version appeared later in the 1970s in papers by Segal [S1, S2] and McDuff [M].  A more immediate predecessor is the paper of Madsen-Tillmann [MT] which explicitly conjectures the Madsen-Weiss theorem.  

The homology of the mapping class group of a closed orientable surface $S$ is the homology of an Eilenberg-MacLane space $K(\Gamma,1)$, with $\Gamma$ the mapping class group.  We are free to choose any $K(\Gamma,1)$ we like, and there is a particular choice that works well for the Madsen-Weiss theorem.  This is the space $\C(S,\R^\infty)$ of all smooth oriented subsurfaces of $\R^\infty$ diffeomorphic to $S$.  The symbol $\C$ is chosen to indicate that $\C(S,\R^\infty)$ is the space of all possible ``configurations" of $S$ in $\R^\infty$.  The space $\C(S,\R^\infty)$ is the union of its subspaces $\C(S,\R^n)$ of smooth oriented subsurfaces of $\R^n$ diffeomorphic to $S$, for finite values of $n$, with the direct limit topology.  Each $\C(S,\R^n)$ is given the usual $C^\infty$ topology as the orbit space of the embedding space $\E(S,\R^n)$ under the action of $\Diff^+(S)$ by composition, where $\Diff^+(S)$ is the group of orientation-preserving diffeomorphisms of $S$ with the $C^\infty$ topology.  By definition, in the direct limit topology on $\C(S,\R^\infty)$ a set is open if and only if it intersects each subspace $\C(S,\R^n) $ in an open set.  A key feature of direct limit topologies (for Hausdorff spaces) is that compact subspaces lie in finite stages of the direct limit, so the homotopy groups and homology groups of a direct limit are the direct limits of the homotopy and homology groups of the finite stages. 

As we explain in Appendix A, the quotient map $\E(S,\R^\infty)\to \C(S,\R^\infty)$ is a fiber bundle with fiber $\Diff^+(S)$, whose total space $\E(S,\R^\infty)$ is contractible, so the long exact sequence of homotopy groups for the bundle gives isomorphisms $\pi_i\Diff^+(S)\iso\pi_{i+1}\C(S,\R^\infty)$ for all $i$.  By the Earle-Eells theorem these groups are trivial for $i>0$, at least when the genus of $S$ is at least $2$, so we see that $\C(S,\R^\infty)$ is a $K(\Gamma,1)$ for the mapping class group in these cases.

To place things in their natural setting, consider a smooth closed orientable manifold $M$ of arbitrary dimension $d\ge 0$, and let  $\C(M,\R^\infty)$ be the space of all smooth oriented submanifolds of $\R^\infty$ diffeomorphic to $M$, defined analogously to the case $d=2$ that we have just described.  
It is natural to ask what information can be extracted from an embedded submanifold $M\subset \R^n$ of dimension $d$ just by looking locally.  Imagine taking a powerful magnifying lens and moving it all around $\R^n$ to see what the submanifold $M$ looks like.  For most positions of the lens one is not close enough to $M$ to see anything of $M$ at all, so one sees just the empty set, but as one moves near $M$ one sees a small piece of $M$ that appears to be almost flat.  Moving toward $M$, this piece first appears at the edge of the lens, then moves to the center.  The space of almost flat $d$-planes in an $n$-ball has the same homotopy type as the subspace of actually flat $d$-planes since one can canonically deform almost flat planes to their tangent planes at their centers of mass.  Regarding the $n$-ball as $\R^n$, we thus have, for each position of the lens where the view of $M$ is nonempty, a $d$-plane in $\R^n$.  This plane has an orientation determined by the given orientation of $M$. We will use the notation $AG_{n,d}$ for the {\it affine Grassmannian\/} of oriented flat $d$-planes in $\R^n$, where the word ``affine" indicates that the planes need not pass through the origin.  Taking into account positions of the lens where the view of $M$ is empty, we then have a point in the one-point compactification $AG^+_{n,d}$ of $AG_{n,d}$ for each position of the lens, with the compactification point corresponding to the empty $d$-plane.  

Positions of the lens near infinity in $\R^n$ give an empty view of $M$, hence map to the compactification point at infinity in $AG^+_{n,d}$.  Thus by letting the position of the lens vary throughout all of $\R^n$ we obtain map $S^n=\R^n \cup \{\infty\}\to AG^+_{n,d}$ taking the basepoint $\infty$ in $S^n$ to the basepoint of $AG^+_{n,d}$, the compactification point.  Such a map is exactly a point in the $n$-fold loopspace $\Omega^n AG^+_{n,d}$.  This point in $\Omega^n AG^+_{n,d}$ associated to the submanifold $M\subset\R^n$ depends on choosing a sufficiently large power of magnification for the lens, which in turn can depend on the embedding of $M$ in $\R^n$.  Making the plausible assumption that the magnification can be chosen to vary continuously with $M$, we then obtain a map $\C(M,\R^n)\to \Omega^n AG^+_{n,d}$.  Letting $n$ increase, the natural inclusion $\C(M,\R^n) \incl \C(M,\R^{n+1})$ corresponds to the inclusion $\Omega^n AG^+_{n,d} \incl \Omega^{n+1}AG^+_{n+1,d}$ obtained by applying $\Omega^n$ to the inclusion $AG^+_{n,d} \incl \Omega AG^+_{n+1,d}$ that translates a $d$-plane in $\R^n$ from $-\infty$ to $+\infty$ in the $(n+1)$st coordinate of $\R^{n+1}$.  Passing to the limit over $n$, we get a map
$$
\C(M,\R^\infty)\to \Omega^\infty AG^+_{\infty,d}
$$
which we will refer to as the {\it scanning map}.  (In the case $d=0$ when $M$ is a finite set of points, the scanning process is described explicitly in [S2], which seems to be the first place where the term ``scanning" is used in this context.)

The affine Grassmannian $AG_{n,d}$ can be described in terms of the usual Grassmannian $G_{n,d}$ of oriented $d$-planes through the origin in $\R^n$. The projection $AG_{n,d} \to G_{n,d}$ translating each $d$-plane to the parallel plane through the origin is a vector bundle whose fiber over a given $d$-plane $P$ is the vector space of vectors orthogonal to $P$ since there is a unique such vector translating $P$ to any given plane parallel to $P$.  The vector bundle $AG_{n,d}$ is thus the orthogonal complement of the canonical bundle over $G_{n,d}$, and the one-point compactification $AG^+_{n,d}$ is the Thom space of this complementary vector bundle since $G_{n,d}$ is compact.  

For the scanning map $\C(M,\R^\infty)\to \Omega^\infty AG^+_{\infty,d}$ the source space depends on the manifold $M$ but the target space does not, so one would hardly expect this map to be any sort of equivalence for arbitrary $M$.  One might have a better chance if one could replace the source by some sort of amalgam or limit over all choices of $M$, and this is what the Madsen-Weiss theorem does in the case of surfaces, when $d=2$, so that there is just one closed orientable surface $S_g$ for each genus $g$.  
We can form a  limit of this sequence of closed surfaces by considering the simplest infinite-genus surface $S_\infty$. This is the union of an increasing sequence of compact surfaces $S_{g,1}$ of genus $g$ with one boundary component. 
\medskip
\vskip4pt
\centerline{\epsfbox{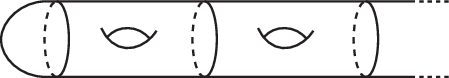}}

\medskip\noindent
Let $\Diff(S_{g,1})$ be the group of diffeomorphisms of $S_{g,1}$ that restrict to the identity on the boundary circle (and whose derivatives of all orders equal those of the identity diffeomorphism at points in the boundary).  The standard embedding $S_{g,1}\incl S_{g+1,1}$ induces a natural inclusion $\Diff(S_{g,1})\incl\Diff(S_{g+1,1})$ by extending diffeomorphisms by the identity on the complement of $S_{g,1}$ in $S_{g+1,1}$.  
Let $\C_g$ be the space of subsurfaces of $(-\infty,g] \times\R^\infty$ diffeomorphic to $S_{g,1}$ whose boundaries coincide (to infinite order) with the boundary of a fixed standard $S_{g,1}\subset (-\infty,g]\times \R^\infty$ with $S_{g,1}\cap (g\times \R^\infty) = \bdy S_{g,1}$.  As in the closed case, $\C_g$ is a $K(\Gamma,1)$ for $\Gamma=\pi_0\Diff(S_{g,1})$.  There are inclusions $\C_g \subset \C_{g+1}$ by adjoining the missing part of the standard $S_{g+1,1}\subset (-\infty,g+1]\times \R^\infty$, so we can form the direct limit $\C_\infty =\cup_g \C_g$.
The version of the Madsen-Weiss theorem that we will prove is the following:

\proclaim The Madsen-Weiss Theorem. There is an isomorphism 
$$
H_*(\C_\infty) \iso H_*(\Omega_0^\infty AG^+_{\infty,2})
$$ 
where $\Omega_0^\infty AG^+_{\infty,2}$ denotes the basepoint path-component of $\Omega^\infty AG^+_{\infty,2}$.

The reason for replacing $\Omega^\infty$ by $\Omega_0^\infty$ is that $\C_\infty$ is path-connected since each $\C_g$ is path-connected.  All the path-components of $\Omega^\infty AG^+_{\infty,2}$ have the same homotopy type since it is an H-space with $\pi_0$ a group, so it does not matter which path-component we choose.

This version of the Madsen-Weiss theorem implies the more usual version in terms of mapping class groups since 
$$
H_i(\C_\infty)=\lim_g H_i(\C_g) = \lim_g H_i(\pi_0\Diff(S_{g,1}))
$$ 
and by Harer stability these limits are actually achieved at finite stages and agree with the homology of mapping class groups of the closed surfaces $S_g$ in the stable dimension range.

One might wish to have a scanning map $\C_\infty \to \Omega^\infty AG^+_{\infty,2}$ that induces the homology isomorphism in the theorem, but constructing a scanning map like this involving surfaces with fixed boundary would require an additional argument and this is not actually needed for the proof of the theorem.  Instead the isomorphism in the theorem will be obtained in a somewhat less direct fashion by recasting the scanning idea in a slightly different form. 

\bigskip\noindent
{\bf A Convention:} \enskip We will sometimes say a map is a homotopy equivalence but only prove it is a weak homotopy equivalence, inducing isomorphisms on all homotopy groups.  This will not be a problem since in the end we are only interested in homology groups and weak homotopy equivalences induce isomorphisms on homology.  Some of the spaces we consider are known to have the homotopy type of CW complexes, when the qualifier ``weak" can be dropped, but this may not be true in all cases.

\section 2. The Proof, Assuming Two Delooping Propositions

Let $\C(M,\R^\infty)$ be the space of oriented smooth submanifolds of $\R^\infty$ diffeomorphic to a given closed orientable manifold $M$ of dimension $d\ge 0$.  By definition, $\C(M,\R^\infty)$ is the union of its subspaces $\C(M,\R^n)$ of oriented smooth submanifolds of $\R^n$ diffeomorphic to $M$, with the direct limit topology.

We enlarge $\C(M,\R^n)$ to the space $\Cn$ of all smooth oriented $d$-dimensional submanifolds of $\R^n$ that are not necessarily closed, but are properly embedded, so they are closed subspaces of $\R^n$ even if they are not closed manifolds. The manifolds in $\C^n$ need not be connected, and the empty manifold is allowed and will be important, serving as a basepoint. The topology on $\Cn$ is defined to have as basis the sets $C(B,V)$ where:
 
\smallskip
\item{(i)} $B$ is a closed ball in $\R^n$ centered at the origin.
 
\item{(ii)} $V$ is an open set in the standard $C^\infty$ topology on the space of $d$-dimensional smooth compact submanifolds of $B$ that are properly embedded (meaning that the intersection of the submanifold with $\bdy B$ is the boundary of the submanifold, and this is a transverse intersection).
 
\item{(iii)} $C(B,V)$ consists of all manifolds $M$ in $\Cn$ meeting $\bdy B$ transversely such that $M\cap B$ is in $V$.  

\smallskip\noindent
The empty submanifold of $B$ is allowed in (ii), and forms an open set $V$ all by itself.  It is not hard to check that the sets $C(B,V)$, as $B$ and $V$ vary, satisfy the conditions for defining a basis for a topology. In this topology a neighborhood basis for a given manifold $M\in\Cn$ consists of the subsets of $\Cn$ obtained by choosing a ball $B$ with $\bdy B$ transverse to $M$ and taking all manifolds $M'\in \Cn$ such that $M'\cap B$ is $C^\infty$-close to $M\cap B$, so in particular $M'$ is also transverse to $\bdy B$.

With this topology, $\Cn$ is path-connected when $n>d$ since a given manifold $M \subset \R^n$ can be connected to the empty manifold by the path obtained by larger and larger radial expansion of $\R^n$ from any point in the complement of $M$.  This path is continuous since in any ball $B$ it is eventually constant, the empty manifold.  If $n=d$, $\Cn$ consists of two points, the empty manifold and $\R^n$ itself. This case is not of interest for us, so {\it we will assume $n>d$ from now on without further mention}.

For the proof of the Madsen-Weiss theorem, when we take $d=2$, it would be possible to restrict surfaces in $\Cn$ to have finite total genus over all components since the constructions we make will preserve this property, but there is no special advantage to adding this restriction.

The homotopy type of $\Cn$ is easy to determine completely:

\proclaim Proposition 2.1. $\Cn$ has the homotopy type of its subspace $AG^+_{n,d}$ consisting of oriented affine linear $d$-planes in $\R^n$, together with the empty set which gives the one-point compactification of the space of $d$-planes.

\pf Proof. Consider the operation of rescaling $\R^n$ from the origin by a multiplicative factor $\lambda\ge 1$, including the limiting value $\lambda=\infty$. Applied to any manifold $M$ defining a point in $\Cn$ this rescaling operation defines a path in $\Cn$ ending with an $d$-plane through the origin if $0\in M$, or with the empty set if $0\notin M$.  This path does not depend continuously on $M$, however, since as we perturb a manifold containing $0$ to one which does not, the endpoint of the path suddenly changes from a plane through the origin to the empty set.  To correct for this problem we will modify the rescaling operation. 

Let $W$ be a tubular neighborhood of a given $M\in\Cn$, so $W$ is a vector bundle $p\:W\to M$ with fibers orthogonal to $M$.  If $0\notin W$ then we do just what we did in the previous paragraph, rescaling $\R^n$ from the origin by factors ranging from $1$ to $\infty$.  If $0\in W$ then we rescale by different factors in the directions of the tangent and normal planes to $M$ at $p(0)$. In the tangential direction we again rescale by a factor going from $1$ to $\infty$, but in the normal direction we rescale by a factor ranging from $1$ to a number $\lambda$ that decreases from $\infty$ to $1$ depending on the position of $0$ in the fiber of $W$, with $\lambda$ being $\infty$ when $0$ is near the frontier of $W$ and $1$ when $0$ is near the zero-section $M$ in $W$, so there is no normal rescaling at all near $M$.  This gives a continuous deformation of $\Cn$ into the subspace $AG^+_{n,d}$ consisting of linear planes and the empty set, using the fact that the tubular neighborhood $W$ can be chosen to vary continously with $M$.  The subspace $AG^+_{n,d}$ is taken to itself during the deformation, so the inclusion $AG^+_{n,d}\incl \Cn$ is a homotopy equivalence.  \qed

We filter $\Cn$ by subspaces $\Cnzero\subset \Cnone\subset\cdots\subset\Cnn=\Cn$ where $\Cnk$ is the subspace of $\Cn$ consisting of manifolds that are properly embedded in $\R^n$ but actually lie in the subspace $\R^k\times (0,1)^{n-k}$.  Thus a manifold $M$ in $\Cnk$ can extend to infinity in only $k$ directions, and the projection $M\to \R^k$ is a proper map, meaning that the inverse image of each compact set in $\R^k$ is compact in $M$.

There is a natural map $\Cnk\to\Omega\Cnkplusone$ obtained by translating a manifold in $\Cnk$ from $-\infty$ to $+\infty$ in the $(k+1)$st coordinate, which gives a loop in $\Cnkplusone$ based at the empty manifold. 
Combining these maps and their iterated loopings gives a composition
$$
\Cnzero \to \Omega \Cnone\to \Omega^2\Cntwo \to\cdots \to \Omega^n\Cn 
$$
This composition takes a closed manifold $M$ in $\Cnzero$ and translates it to infinity in all directions.  This is essentially the same as scanning $M$ with an infinitely large lens.  By shrinking the lens to a small size one obtains a homotopy from this map $\Cnzero\to \Omega\Cn$ to the original scanning map.  However, this observation will not play a significant role in the proof of the Madsen-Weiss theorem. Instead, the proof will consist of two main steps.  The first one is easier and is valid for all $d\ge 0$:

\smallskip
\item{(1)} The map $\Cnk\to\Omega\Cnkplusone$ is a homotopy equivalence when $k > 0$. 
\smallskip
\noindent
This leaves only the map $\Cnzero\to\Omega\Cnone$ to be considered.  This is considerably more delicate except in the one easy case $d=0$ which we discuss in Section~5.  A rough statement of the second step in the case of the Madsen-Weiss theorem is:

\smallskip
\item{(2)} When $d=2$ the map $\Cnzero\to\Omega\Cnone$ becomes a homology equivalence after passing to suitable limits involving letting $n$ and the genus $g$ to go to infinity, so we obtain isomorphisms $H_*(\C_\infty)\iso H_*(\Omega_0 \Cinfone)$.

\medskip

In order to prove (1) and (2) we will not work directly with loopspaces but rather with classifying spaces, more specifically with classifying spaces of topological monoids.  Recall that a topological monoid is a space $\M$ with a continuous product $\M\times\M\to\M$ which is associative and has a two-sided identity element, but need not have inverses.  
The space $\Cnk$ is an H-space when $k<n$, with the product given by juxtaposition in the $(k+1)$st coordinate, after compressing the interval $(0,1)$ in this coordinate to $(0,1/2)$ as in the definition of composition of loops. The homotopy-identity element for the H-space structure is the empty configuration, and the multiplication is homotopy-associative.  Just as the Moore loopspace provides a monoid version of the usual loopspace, with strict associativity and a strict identity, there is a monoid version $\Mnk$ of $\Cnk$. This is the subspace of $\Cn\times [0,\infty)$ consisting of pairs $(M,a)$ with $M$ a manifold in $\R^k\times (0,a)\times (0,1)^{n-k-1}$. The product in $\Mnk$ is again given by juxtaposition in the $(k+1)$st coordinate, but this time without any compression.  It is easy to see that the inclusion $\Cnk\incl\Mnk$ as pairs $(M,1)$ is a homotopy equivalence.

A topological monoid $\M$ has a classifying space $B\M$.  Let us recall the construction.  In the special case that $\M$ has the discrete topology, $B\M$ is the $\Delta$-complex having a single vertex, an edge for each element of $\M$, and more generally a $p$-simplex for each $p$-tuple $(m_1,\cdots,m_p)$ of elements of $\M$.  The faces of such a simplex are the $(p-1)$-tuples obtained by deleting the first or last $m_i$ or by replacing two adjacent $m_i$'s by their product in $\M$.  Thus $B\M$ is a quotient space of $\coprod_p \Delta^p \times \M^p$ with certain identifications over $\bdy\Delta^p \times \M^p$ for each $p$.  Essentially the same construction can be made when $\M$ has a nontrivial topology.  One gives each $ \Delta^p \times \M^p$ the product topology and then one forms a quotient of $\coprod_p \Delta^p \times \M^p$ using the same rules for identifications over the subspaces $\bdy\Delta^p \times \M^p$ as before.  Thus $B\M$ is built from $p$-simplices for all $p\ge 0$, where the set of $p$-simplices for fixed $p$ is not a discrete set but forms a space $\M^p$.  In particular there is a single $0$-simplex.

There is a natural map $\M\to\Omega B\M$ obtained by taking $p=1$ in the definition of $B\M$, so each element  $m\in\M$ gives an edge $\Delta^1\times m$ in $B\M$ which is a loop at the basepoint vertex of $B\M$.  In favorable cases this map $\M\to\Omega B\M$ is a homotopy equivalence.  This happens when $\pi_0\M=0$ for example, and more generally when $\pi_0\M$ is a group with respect to the multiplication coming from the monoid structure.  A proof of this classical fact is given in Appendix D.  

In Section 3 we will prove the following result which will complete step (1) of the proof:

\proclaim Proposition 3.1.  For $k>0$ and $d\ge 0 $ the monoid $\pi_0\Mnk$ is a group and $B\Mnk \simeq\Cnkplusone_0$, the component of $\Cnkplusone$ containing the empty manifold. Hence $\Cnk\simeq \Omega\Cnkplusone$ when $k>0$ and $d\ge 0$.

Here the last statement follows from the first two by the chain of equivalences 
$$
\Cnk\simeq \Mnk \simeq \Omega B\Mnk \simeq \Omega \Cnkplusone_0 = \Omega\Cnkplusone
$$
where the last equality holds since $\Omega X$ only involves the basepoint component of $X$.  The proof of Proposition 3.1 will show in fact that the natural map $\Cnk\to\Omega\Cnkplusone$ is a homotopy equivalence.

In the more delicate case $k=0$ the monoid $\pi_0\Mnzero$ is not a group since all the manifolds in each path component of $\Mnzero$ are diffeomorphic and the disjoint union of two nonempty manifolds is not empty.  For an arbitrary topological monoid $\M$, if the map $\M\to\Omega B\M$ is a homotopy equivalence then $\pi_0\M$ must be a group since the induced map $\pi_0\M\to\pi_0\Omega B\M$ is a homomorphism and $\pi_0\Omega B\M =\pi_1 B\M$ is a group.  Thus the map $\Cnzero\to \Omega\Cnone$ cannot be a homotopy equivalence.  

Another issue is that the product operation in $\Mnzero$ is disjoint union whereas what one needs in the Madsen-Weiss theorem is a sum operation within the realm of connected surfaces.  With this in mind we will replace $\Mnzero$ in the case $d=2$ by a monoid of connected surfaces with nonempty boundary, so that gluing two surfaces together by identifying a boundary circle of one with a boundary circle of the other preserves the connectedness property.

In order to define the new monoid, we first consider the cylinder $Z = \R\times C \subset \R\times \R^{n-1}$ where $C$ is a fixed circle in the first two coordinates of $\R^{n-1}$ centered at the origin and of small radius, say radius $1/2$.  Then we let $\Mn$ be the space of pairs $(S,a)$ where $S$ is a compact connected orientable surface in $[0,a]\times(-1,1)^{n-1} $ such that:

\smallskip
\item {(i)} $\bdy S$ consists of the two circles $Z\cap \bigl(\{0,a\} \times (-1,1)^{n-1}\bigr)$, and $S$ is tangent to $Z$ to infinite order along these two circles. Furthermore, $S-\bdy S \subset (0,a)\times (-1,1)^{n-1}$.
\item{(ii)} $S\cap \bigl([0,a] \times (-1,0]\times (-1,1)^{n-2}\bigr) = Z\cap \bigl([0,a] \times (-1,0]\times (-1,1)^{n-2}\bigr)$. 

\vskip8pt
\centerline{\epsfbox{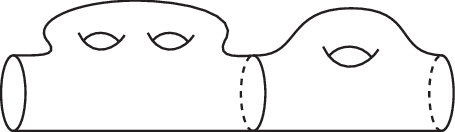}}
\vskip-2pt
\smallskip\noindent
We also allow the limiting case $a=0$ when $S$ degenerates to just a circle. The monoid structure in $\Mn$ is given by juxtaposition in the first coordinate of $\R^n$.  The condition of tangency to infinite order in (i) guarantees that the monoid operation stays within the realm of smooth surfaces.

We will only need the case $n=\infty$, so we set $\Minf=\cup_n\Mn$ with the direct limit topology as usual.  In Section~4 we will prove:

\proclaim Proposition 4.1.  $B\Minf\simeq\Cinfone $.

The path-components of $\Minf$ are the subspaces $\Minf_g$ consisting of surfaces of genus $g$, so $\pi_0\Minf$ is the monoid $\Z_{\ge 0}$ of nonnegative integers under addition.  We have seen that for a monoid $\M$ such that $\pi_0\M$ is not a group, the map $\M\to \Omega B\M$ cannot be a homotopy equivalence since $\pi_0\Omega B\M$ is always a group, namely $\pi_1 B\M$.  It is not hard to see that $\pi_1 B\M$ is the group completion of $\pi_0\M$, which can be defined as the group with a presentation consisting of the elements of $\pi_0\M$ as generators and the multiplication table of $\pi_0\M$ as relations.  For example if $\pi_0\M$ is the monoid $\Z_{\ge 0}$ of nonnegative integers under addition then its group completion is just $\Z$.

When $\pi_0\M$ is only a monoid there is a theorem, known as the Group Completion Theorem, that describes the homology of $\Omega B\M$ in terms of the homology of $\M$, under certain hypotheses on $\M$.  To state this let us assume for simplicity that $\pi_0\M$ is $\Z_{\ge 0}$.  The path-components $\M_p$ of $\M$ then correspond to integers $p\ge 0$ and there is a stabilization map $\M_p \to \M_{p+1}$ given by taking the product (on the right, say) with any element of $\M_1$. In our application this map will be an injection $\M_p \incl \M_{p+1}$ so one can form the direct limit of the $\M_p$'s as $\cup_p \M_p$, with the direct limit topology. The Group Completion Theorem then says that if the multiplication in $\M$ is homotopy-commutative, then there is an isomorphism 
$$
H_*(\Z\times\cup_p \M_p)\iso H_*(\Omega B\M)
$$
On the left side we have $\pi_0(\Z\times\cup_p \M_p) = \Z $ since $\cup_p\M_p$ is path-connected, as each $\M_p$ is path-connected.  Thus passing from $\M$ to $\Z \times\cup_p \M_p$ has the effect on $\pi_0$ of converting $\Z_{\ge 0}$ to its group completion $\Z$.  It also has the effect of forcing all the path-components of $\Z \times\cup_p \M_p$ to be homotopy equivalent, just as they are in $\Omega B\M$ since a loopspace is an H-space with $\pi_0$ a group.

We can apply the Group Completion Theorem to $\M=\Minf$ since it is easy to see that the multiplication in $\Minf$ is homotopy-commutative, using the same sort of idea as is used to show that higher homotopy groups are abelian after $\pi_1$.  Restricting to a single path-component, we then obtain isomorphisms $H_*(\cup_g\Minf_g)\iso H_*(\Omega_0 B\Minf)$.

The isomorphism in the Madsen-Weiss Theorem is now obtained by composing seven isomorphisms:
$$\eqalignno{
H_*(\C_\infty) &\iso H_*(\cup_g \Minf_g) &(1)\cr
&\iso H_*(\Omega_0 B\Minf) &(2)\cr
&\iso H_*(\Omega_0\Cinfone) &(3)\cr
&\iso \lim_n H_*(\Omega_0\Cnone) &(4)\cr
&\iso \lim_n H_*(\Omega^n_0\Cn) &(5)\cr
&\iso \lim_n H_*(\Omega^n_0 AG^+_{n,2}) &(6)\cr
&\iso H_*(\Omega^\infty_0 AG^+_{\infty,2}) &(7) \cr
}$$
The first isomorphism follows from the fact that $\C_\infty =\cup_g \C_g$ and $\C_g \simeq \Minf_g$ (and in fact $\C_g$ and $\Minf_g$ are homeomorphic).  
The second isomorphism comes from the Group Completion Theorem as we have just noted.  
The third isomorphism follows from Proposition~4.1. 
The fourth isomorphism holds since homology and the loopspace functor commute with direct limit for a space that is the union of a sequence of subspaces with the direct limit topology.  The fifth isomorphism follows by repeated applications of Proposition~3.1.  The sixth is Proposition~2.1, and the seventh follows in the same way as the fourth.

\medskip

\section 3. Delooping: the Easy Cases

In this section we prove the first of the two missing steps from the proof in Section~2:

\proclaim Proposition 3.1.  For $k>0$ and $d\ge 0 $ the monoid $\pi_0\Mnk$ is a group and $B\Mnk \simeq\Cnkplusone_0$, the component of $\Cnkplusone$ containing the empty manifold. Hence $\Cnk\simeq \Omega\Cnkplusone$ when $k>0$ and $d\ge 0$.

First we give a more precise version of the first assertion.

\proclaim Lemma 3.2. $\pi_0\Cnk = 0$ when $k > d $, while for $0<k\le d$ $\pi_0\Cnk$ is isomorphic to the group $\Omega^{SO}_{d-k,n-k}$ of cobordism classes of closed oriented $(d-k)$-manifolds in $\R^{n-k}$, where cobordisms are embedded in $\R^{n-k} \times I$.

The sum operation in $\Omega^{SO}_{d-k,n-k}$ is given by disjoint union, after translating the two submanifolds of $\R^{n-k}$ to lie on opposite sides of a hyperplane.  The inverse of a given manifold $M\subset \R^{n-k}$ is obtained by reflecting across a hyperplane, where a cobordism between the union of $M$ and its reflection is given by a U-shaped embedding of $M\times I $ in $\R^{n-k} \times I$, with both ends of $M\times I$ lying in the same end of $\R^{n-k} \times I$.

For the Madsen-Weiss theorem we only need these cobordism groups when $d=2$ and $k=1,2$.  The actual calculation of these groups will not be needed, but it is not hard to see that they are given by $\Omega^{SO}_{0,n-2} = \Z$ and $\Omega^{SO}_{1,n-1}=0$.

\medskip
\pf Proof of 3.2. A point in $\Cnk$ is a $d$-dimensional manifold $M$ in $\R^n$ with the projection $p\: M\to\R^k$ a proper map. By a small perturbation of $M$ we can arrange that $p$ is transverse to $0\in\R^k$. If $k>d$ this implies that $p^{-1}(0)$ is empty.  Then there is a ball $B$ about $0$ in $\R^k$ disjoint from $p(M)$ since $p(M)$ is a closed subset of $\R^k$, the map $p$ being proper.  Expanding the $\R^k$ factor of $\R^n$ radially from $0$ until $B$ covers all of $\R^k$ then gives a path in $\Cnk$ from $M$ to the empty manifold, so $\pi_0\Cnk=0$ in the cases $k>d$.

When $k \le d$ transversality and properness of $p$ imply that $p^{-1}(0)=M\cap \bigl(\{0\}\times\R^{n-k}\bigr)$ is a closed submanifold $M_0\subset M$ of dimension $d-k$.  The manifold $M_0$ inherits an orientation from the given orientation of $M$ and an orientation of $\R^k$ which pulls back to an orientation of the normal bundle of $M_0$ in $M$.  Thus $M_0$ gives an element of $\Omega^{SO}_{d-k,n-k}$.  By a standard transversality argument, the association $M\mapsto M_0$ determines a well-defined map $\varphi\: \pi_0\Cnk\to \Omega^{SO}_{d-k,n-k}$, and this is a homomorphism since the sum is disjoint union in both cases. Moreover, $\varphi$ is surjective since we can choose $M=\R^k\times M_0$ for any cobordism class $[M_0]\in \Omega^{SO}_{d-k,n-k}$.

To show that $\varphi$ is injective, first deform a given $M$ near $M_0$ so that $M$ agrees with $\R^k\times M_0$ in a small neighborhood of $M_0$. Then by radial expansion of $\R^k$ we can construct a path in $\Cnk$ connecting $M$ to $\R^k\times M_0$, so we may assume $M=\R^k\times M_0$.  If we are given another $M'\in\Cnk$ with $[M_0]=[M'_0]$ in $\Omega^{SO}_{d-k,n-k}$, we may similarly assume $M'=\R^k \times M'_0$. Since $[M_0]=[M'_0]$ there is a cobordism $V\subset I \times \R^{n-k}$ from $M_0$ to $M'_0$, and we can assume this lies in $I\times (0,1)^{n-k}\subset I \times \R^{n-k}$.  Let $W$ be the submanifold $\R^{k-1}\times V \subset\R^n$, extended by $M$ in $\R^{k-1}\times (-\infty,0]\times \R^{n-k}$ and by $M'$ in $\R^{k-1}\times [1,\infty)\times \R^{n-k}$. Translating $W$ from $+\infty$ to $-\infty$ in the $k$th coordinate of $\R^n$ then gives a path in $\Cnk$ from $M$ to $M'$. This shows $\varphi$ is injective.  \qed

\centerline{\epsfbox{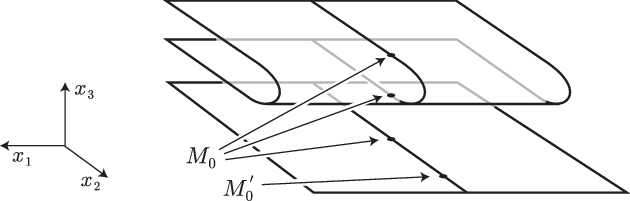}}

Now we turn to proving the second assertion of Proposition 3.1, the homotopy equivalence $B\Mnk\simeq \Cnkplusone_0$.  We will define a map $\sigma\: B\Mnk\to \Cnkplusone_0$ and show this induces isomorphisms on all homotopy groups.  A point in $B\Mnk$ is given by a $p$-tuple of points $m_1,\cdots,m_p$ in $\Mnk$ and a point in $\Delta^p$ with barycentric cooordinates $w_0,\cdots,w_p$.  When we form the product $m_1\cdots m_p$ we obtain manifolds $M_1,\cdots,M_p$ whose $(k+1)$st coordinates lie in intervals $[a_0,a_1],[a_1,a_2],\cdots,[a_{p-1},a_p]$ for $0=a_0\le a_1\le\cdots\le a_p$. The union of the manifolds $M_i$ is a manifold $M\in\Cnkplusone$.  However, $M$ does not depend continuously on the given point of $B\Mnk$ since as the weight $w_0$ or $w_p$ goes to $0$, the manifold $M_1$ or $M_p$ is suddenly deleted from $M$.  Furthermore, when $w_0$ goes to $0$ the remaining manifold is suddenly translated a distance $a_1$ to the left in the $(k+1)$st coordinate of $\R^n$.  The latter problem can easily be resolved by dropping the restriction $a_0=0$ and translating in the $(k+1)$st coordinate so that the barycenter $b =\sum_i w_i a_i$ is at $0$. 

To fix the other problem we choose ``upper and lower barycenters" $b^+$ and $b^-$ by letting $a^+_i=\max\{a_i,b\}$ and $a^-_i=\min\{a_i,b\}$, then setting $b^+=\sum_i w_i a^+_i$ and $b^- = \sum_i w_i a^-_i$.  Then we have $a_0\le b^-\le b \le b^+ \le a_p$, with all these inequalities strict unless they are all equalities (and hence $M$ is empty).  
Now we define the map $\sigma\: B\Mnk\to \Cnkplusone_0$ by taking the part of $M$ in the open slab $S(b^-,b^+)=\R^k\times (b^-,b^+)\times \R^{n-k-1}$ and stretching this slab out to $\R^n$ by stretching the interval $(b^-,b^+)$ to $(-\infty,+\infty)$ in the $(k+1)$st coordinate. More precisely, let $L\:\R^n \to \R^n$ be the map that is the identity in all coordinates except the $(k+1)$st coordinate, where it is the composition of a fixed identification of $(-\infty,+\infty)$ with $(-1,+1)$ followed by the unique affine linear map taking $[-1,+1]$ onto $[b^-,b^+]$.  Then $\sigma$ of the given point in $B\Mnk$ is defined to be $L^{-1}(M)$.  Note that this definition works even in the degenerate case that $b^+=b^-=0$ when we are just stretching the empty manifold.  The map $\sigma$ is continuous since it is continuous on each product $\Delta^p\times (\Mnk)^p$ and is consistently defined when weights go to $0$.  The image of $\sigma$ lies in the path-component $\Cnkplusone_0$ of $\Cnkplusone$ containing the empty manifold since the image contains the empty manifold and $B\Mnk$ is path-connected.

To show that $\sigma_*\:\pi_q(B\Mnk)\to\pi_q(\Cnkplusone)$ is surjective, represent an element of the target group by a map $f\:D^q\to \Cnkplusone$ taking $\bdy D^q$ to the basepoint of $\Cnkplusone$, the empty manifold. Thus we have manifolds $M_t = f(t)\in \Cnkplusone$ for $t\in D^q$.  For the moment let us fix a value of $t$ and let $M=M_t$.  
The projection $p\:M \to \R^{k+1}$ is a proper map, so its set of regular values is open and dense in $\R^{k+1}$.  Let $x$ be a regular value and let $M_x=p^{-1}(x)$, a manifold of dimension $d-k-1$.  
In the cases $d<k+1$ the manifold $M_x$ is empty.  Then we can choose a closed ball $B\subset \R^k$ and a closed interval $J\subset\R$ such that $B\times J$ is disjoint from $p(M)$ in $\R^{k+1}$.  By expanding each ball $B\times \{a\}$ for $a\in J$ to $\R^k\times \{a\}$, damping this expansion down to zero near $\bdy J$, we obtain a deformation of $M$ to a manifold disjoint from the slices $S(a)=\R^k\times \{a\}\times (0,1)^{n-k-1}$ for $a$ in a subinterval of $J$.  

When $d\ge k+1$ the manifold $M_x$ can be nonempty, and we will achieve the same disjointness property by using a process somewhat like the one used in the proof of Lemma~3.2.  First perturb $M$  to make it agree with $\R^{k+1}\times M_x$ over a neighborhood $B\times J$ of $x$ .  Next, do the damped radial expansion in slices $S(a)$ as in the cases $d<k+1$ to make $M$ agree with $\R^{k+1}\times M_x$ in all slices $S(a)$ for $a$ in a subinterval $J'$ of $J$.  Since $M$ is in the component $\Cnkplusone_0$ of $\Cnkplusone$, the cobordism class $[M_x]$ is zero in $\Omega^{SO}_{d-k-1,n-k-1}$, so $M_x$ bounds a manifold $V$ in $I\times \R^{n-k-1}$. We place the $I$ factor of this product in the $k$th coordinate of $\R^k$.  Then translation in this coordinate as in the proof of Lemma~3.2, but now damped down near $\bdy J'$, gives a further deformation of $M$ to make it disjoint from all slices over a subinterval of $J'$.

\vskip8pt
\centerline{\epsfbox{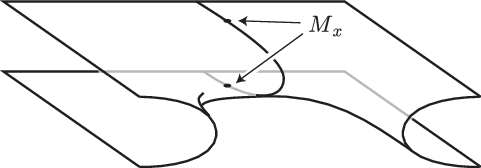}}
\vskip-2pt
For nearby values of $t$ we can use the same product $B\times J$.  Using compactness of $D^q$ we can then choose a cover of $D^q$ by finitely many open sets $V_i$ with corresponding products $B_i\times J_i$ as above for $t\in V_i$.  After shrinking the intervals $J_i$ appropriately, we can assume they are all disjoint.  Then we can perform the deformations of $M_t$ for different products $B_i\times J_i$ independently, damping the deformation for $B_i\times J_i$ down to zero outside slightly smaller open sets $U_i$ that still cover $D^q$.  Performing these damped down deformations for all the products $B_i\times J_i$ produces a new map $f\: D^q \to \Cnkplusone$ homotopic to the old one by a homotopy which is constant for $t\in\bdy D^q$ since the manifold $M_t$ is empty there.  We discard the original $f$ and work now with the new family $f(t)=M_t$.

We can choose slices $S(a_i)$ disjoint from $M_t$ over $U_i$ with $a_i\in J_i$, and we can choose corresponding weights $w_i$ via a partition of unity subordinate to the cover $\{U_i\}$ of $D^q$.  Over a neighborhood of $\bdy D^q$ we can choose just the slice $S(0)$ with weight $1$.  By taking the parts of $M_t$ in the slabs between adjacent slices $S(a_i)$ we obtain a map $g\:D^q \to B\Mnk$ taking $\bdy D^q$ to the basepoint.  The composition $\sigma g$ is homotopic to $f$ by expanding the intervals $(b^-_t,b^+_t)$ for the family $g(t)$ to $(-\infty,+\infty)$, thereby deforming the maps $L_t$ to the identity.

The argument for showing $\sigma$ is injective on $\pi_q$ is similar.  Here we start with a map $g\: S^q\to B\Mnk$ and a map $f\:D^{q+1}\to \Cnkplusone$ restricting to $\sigma g$ on $S^q$.  We deform the family $f(t)=M_t$ by the same procedure as before, expanding in various slices $S(a)$ to make $M_t$ disjoint from slices $S(a_i)$ over sets $U_i$ covering $D^{q+1}$, with weights $w_i$ coming from a partition of unity supported in the cover $\{U_i\}$.  Over $S^q$ this deformation of $f$ induces a corresponding deformation of $g$ since slices are preserved during the deformation.  We use the same notations for the new $f$ and $g$.  The maps $L_t$ that occur in the definition of $\sigma g$ are initially defined over $S^q$, and we extend them so that they are defined over all of $D^{q+1}$ by deforming them to the identity as $t$ moves across a collar neighborhood of $S^q$ in $D^{q+1}$, then taking them to be the identity over the rest of $D^{q+1}$.  We then obtain a map $f' \: D^{q+1}\to B\Mnk$ by taking the parts of $L_t(M_t)$ between the slices $L_t(S(a_i))$, with the weights $w_i$ on these slices.  A homotopy from $f$ to $f'$ on $S^q$ is obtained by letting the weights on the slices for $f$ decrease to $0$ while the weights on the slices $L_t(S(a_i))$ for $f'$ increase from $0$ to their chosen values $w_i$.  Thus $f$ is homotopic to a map that extends over $D^{q+1}$, which shows that $\sigma$ is injective on $\pi_q$.  \qed

\section 4. Delooping: the Harder Case

In this section we restrict to the case of surfaces, so $d=2$.  We also take $\Cinfone$ to consist of surfaces in $\R\times (-1,1)^\infty$ rather than $\R\times (0,1)^\infty$.  This does not change the homeomorphism type of $\Cinfone$, and the new $\Cinfone$ fits better with the definition of $\Minf$, so there is a map $\sigma\: B\Minf\to \Cinfone$ defined just like the $\sigma$ in the preceding section.

\proclaim Proposition 4.1. The map $\sigma\:B\Minf\to\Cinfone $ is a homotopy equivalence.

\pf Proof. We break the argument into two parts by defining a space $\Cinfone_s$ with maps
$$
B\Minf \mapleft\rho \Cinfone_s \mapright\tau \Cinfone
$$
such that the composition $\sigma\rho$ is homotopic to $\tau$, and then we show that $\rho$ and $\tau$ are homotopy equivalences.  This implies that $\sigma$ is also a homotopy equivalence.

A point in $\Cinfone_s$ consists of an oriented surface $S$ in $\Cinfone$ with numbers $a_0\le\cdots\le a_p$, $p\ge 0$, and corresponding weights $w_i\ge 0$ summing to $1$, such that:
\smallskip

\item {(i)} $S\cap \bigl(\R \times (-1,0]\times (-1,1)^\infty\bigr) = Z\cap \bigl(\R \times (-1,0]\times (-1,1)^\infty\bigr)$, a half-cylinder $B$ that we call the {\it base\/} of $S$.  Furthermore, the orientation of $S$ agrees with a fixed orientation of $B$.

\item{(ii)} $S\cap S(a_i)=Z\cap S(a_i)$ for each $i$ and $S$ is tangent to $Z$ to infinite order along these circles. 

\item{(iii)} The intersection of $S$ with each slab between adjacent slices $S(a_i)$ and $S(a_{i+1})$ is connected. 

\smallskip\noindent
When a weight $w_i$ is zero the corresponding $a_i$ can be deleted. 
There is a map $\rho\: \Cinfone_s\to B\Minf$ obtained by restricting to the parts of surfaces between adjacent slices, and there is a map $\tau\: \Cinfone_s\to \Cinfone$ forgetting the extra data. The composition $\sigma\rho$ is homotopic to $\tau$ by expanding the intervals $[b^-,b^+]$ in the definition of $\sigma$ to $[-\infty,+\infty]$.

First we show that $\rho$ induces isomorphisms on all homotopy groups.  For this it suffices to show that for every map $f\:D^q\to B\Minf$ with a lift $\tilde f\: \bdy D^q\to \Cinfone_s$ there is a homotopy of $f$ and a lifted homotopy of $\tilde f$ to new maps $g$ and $\tilde g$ such that $\tilde g$ extends to a lift of $g$ over all of $D^q$, for arbitrary $q> 0$. (Both $B\Minf$ and $\Cinfone_s$ are path-connected so there is no need to consider the case $q=0$.)  

As in the definition of $\sigma$, we view $f$ as defining a family of surfaces $S_t$ in $\R^\infty$ with slices $S(a_i(t))$ and weights $w_i(t)$, for $t\in D^q$.  The lift $\tilde f$ over $\bdy D^q$ gives an enlargement of the surfaces $S_t$ so that they extend to $\pm\infty$ in the first coordinate.  We can assume the barycenters $b(t)$ are always at $0$.  The slices $S(a_i(t))$ cut $S_t$ into pieces, and the first step in the deformation of $f$ and $\tilde f$ is to spread these pieces apart by inserting ``padding" consisting of a piece of the cylinder $Z$ of width $\varepsilon w_i(t)$ at the slice $S(a_i(t))$ for each $i$. This spreads the slices apart, and we take the new $S(a_i(t))$ to be at the center of the inserted piece of $Z$.  Again we translate so that the new barycenters are at $0$.  Letting $\varepsilon$ go from $0$ to $1$ gives a deformation of $f$ and $\tilde f$ to new maps $f$ and $\tilde f$ which we use for the next step of the argument.  

We replace the construction of the intervals $(b^-,b^+)$ by a different construction of intervals $(c^-(t),c^+(t))$ that contain at least one point $a_i(t)$ as follows.  For fixed $t\in D^q$ choose slices $S(a^-(t))$ and $S(a^+(t))$ from among the slices $S(a_i(t))$ with nonzero weights $w_i(t)$, such that $a^-(t)\le 0\le a^+(t)$. Expand the interval $[a^-(t),a^+(t)]$ slightly to an interval $[\phi^-(t),\phi^+(t)]$ by including a little of the padding around $a^-(t)$ and $a^+(t)$. The same choices of $a^\pm(t)$ and $\phi^\pm(t)$ work in a small neighborhood $U$ of $t$-values. Extend the functions $\phi^\pm(t)$ over all of $D^q$ by letting them be $0$ outside a slightly larger neighborhood $V$ and interpolating continuously in $V-U$.  Since $D^q$ is compact, finitely many such neighborhoods $U$ suffice to cover $D^q$. Call these neighborhoods $U_j$, with corresponding intervals $[a_j^-(t),a_j^+(t)]$ and $[\phi_j^-(t),\phi_j^+(t)]$.  Now define the interval $[c^-(t),c^+(t)]$ to be $\cup_j[\phi_j^-(t),\phi_j^+(t)]$.  The key properties of $[c^-(t),c^+(t)]$ are that it is contained in a small neighborhood of $[a_0(t),a_p(t)]$ and it contains $[a_j^-(t),a_j^+(t)]$ for $t\in U_j$.  

Let $g(t)$ be the family obtained from $f(t)$ by replacing the slices $S(a_i(t))$ with the slices $S(a_j^\pm(t))$, with weights $w_j^\pm(t)$ given by a partition of unity for the cover of $D^q$ by the neighborhoods $U_j$.  A homotopy from the family $f(t)$ to the family $g(t)$ is obtained by letting the weights go linearly from their values for $f(t)$ to their values for $g(t)$. Over $\bdy D^q$ this lifts to a homotopy of $\tilde f(t)$ to a lift $\tilde g(t)$ of $g(t)$ since we are not changing $S_t$.  Finally, by expanding the intervals $(c^-(t),c^+(t))$ to $(-\infty,+\infty)$ we obtain homotopies of $g(t)$ and $\tilde g(t)$ to $h(t)$ and $\tilde h(t)$ such that $ h(t)$ lifts to $\tilde h(t)$ over all of $D^q$.  Thus $\rho$ induces isomorphisms on all homotopy groups.

\smallskip
It remains to see that the map $\tau\:\Cinfone_s\to\Cinfone$ is a homotopy equivalence.    Its image is contained in the subspace $\Cinfone_b \subset\Cinfone$ consisting of surfaces satisfying just the `base' condition (i) above.  It is easy to see that the inclusion map $i\: \Cinfone_b\to \Cinfone$ is a homotopy equivalence by constructing a homotopy inverse $j \: \Cinfone\to\Cinfone_b$ as follows.  For $S\in \Cinfone$ let $S'$ be the surface obtained from $S$ by compressing the interval of the second coordinate $x_2$ from $(-1,1)$ to $(1/2,1)$.  Then let $j(S)=S'\cup Z$.  A homotopy from $ij$ to the identity is obtained by first realizing the $x_2$-compression via an isotopy, then introducing the cylinder $Z$ by an isotopy, namely, take the part of $Z$ with $x_1\ge t$ and cap off its boundary circle with a disk (smoothly), then let $t$ go from $+\infty$ to $-\infty$ to obtain the isotopy.  A homotopy from the other composition $ji$ to the identity can be constructed as an isotopy from $S'\cup Z$ to $S$ for $S\in \Cinfone_b$ which `zips' $S'$ to $Z$ by the isotopy shown in the figure below, moving the `zipper' from $+\infty$ to $-\infty$ and re-expanding $(S-B)'$ back to $S-B$.

\vskip12pt
\centerline{\epsfbox{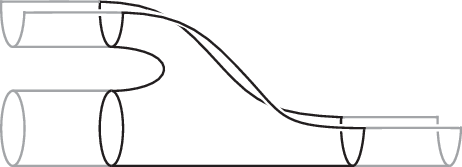}}

To show that the map $\Cinfone_s\to\Cinfone_b$ induces an isomorphism on homotopy groups we start with $f\:D^q\to\Cinfone_b$ and a partial lift $\tilde f\:\bdy D^q\to \Cinfone_s$. There will be three steps.  As notation we let $S[a,b]$ denote the slab between the slices $S(a)$ and $S(b)$.

\noindent
(1) \enskip In this easy preliminary step we deform $f$ and $\tilde f$ so that the slices $S(a_i(t))$ in the family $\tilde f$ all lie in the slab $S[-1,1]$. First we deform $\tilde f$ so that the slices $S(a_i(t))$ vary continuously not just over the set where the corresponding weights $w_i(t)$ are nonzero, but over the closure of this set as well.  We do this by modifying the weights by stretching the interval$[0,\max\{w_i(t)\}]$ to $[-\max\{w_i(t)\},\max\{w_i(t)\}]$, discarding slices whose weights are then negative, and rescaling so that the new weights sum to $1$.  After doing this we do translations in the first coordinate of $\R^\infty$ to move the barycenter of the set of slices $S(a_i(t))$ to $0$ for all $t\in\bdy D^q$.  Then we rescale the first coordinate to squeeze all these slices into $S[-1,1]$, using the fact that the values $a_i(t)$ are bounded as $t$ ranges over the compact space $\bdy D^q$.  These translations and rescalings over $\bdy D^q$ can be extended over $D^q$, damping them off as we move into the interior of $D^q$.

\smallskip\noindent
(2) \enskip Next we deform the family $S_t$ so that for each $t$, the subspace $S_t \cap S[-1,1]$ of $S_t$ is contained in the base component of $S_t$, the component containing the half-cylinder $B$.  Thus, each point of $S_t\cap S[-1,1]$ will be joinable to $B$ by a path in $S_t$, although this path need not stay in $S[-1,1]$.  We can achieve this as follows.  For fixed $t$ choose a slab $S[a,b]$ slightly larger than $S[-1,1]$ with $S(a)$ and $S(b)$ transverse to $S_t$.   
If a component $C$ of $S_t \cap S[a,b]$ does not lie in the base component of $S_t$, we can create a tube in $S[a,b]$ joining $C$ to the base component by the deformation shown in the figure below which takes a pair

\centerline{\epsfbox{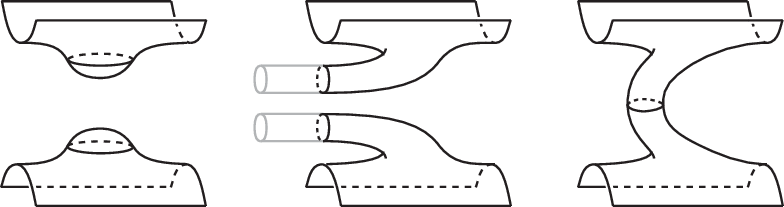}}
\vskip3pt
\noindent
of small disks in $S_t$, one in $C$ and the other near $B$, and drags these disks to $\pm\infty$ in the first coordinate of $\R^\infty$ to create a pair of tubes to $\pm\infty$, then these tubes are brought back from $\pm\infty$ joined together so that the given orientation of $S_t$ extends over the new tube.  
We can do this simultaneously for each such component $C$, and after this is done, all of $S_t \cap S[-1,1]$ lies in the base component of $S_t$ since $S_t \cap S[a,b]$ lies in the base component. 

To do this procedure for all $t$ we first use compactness of $D^q$ to choose a finite cover of $D^q$ by open sets $U_i$ with slabs $S[a_i,b_i]$ containing $S[-1,1]$ such that the slices $S(a_i)$ and $S(b_i)$ are transverse to $S_t$ over the closure $\overline U_i$ of $U_i$. The pairs of disks that we use to create tubes making $S_t\cap S[a_i,b_i]$ connected will be small neighborhoods of pairs of points $p_{ij}$ and $q_{ij}$, with 
$q_{ij}$ near $B$.  The main concern will be choosing all these points to be disjoint for fixed $t$ and varying $i$ and $j$.  This is easy for the points $q_{ij}$ which lie near the fixed $B$, so we focus attention on the points $p_{ij}$.  We choose the $p_{ij}$ by induction on $i$.  For the induction step of choosing $p_{ij}$ over $\overline U_i$, note first that the surface $S_t \cap S[a_i,b_i]$ varies by isotopy as $t$ ranges over $\overline U_i$, so we can regard $S_t \cap S[a_i,b_i]$ as being independent of $t$ over $\overline U_i$.  Inductively we assume we have already chosen points $p_{kl}$ for $k<i$, each point $p_{kl}$ being defined and varying continuously over some closed ball in $\overline U_k$, with all these points $p_{kl}$ being disjoint for each $t$.  To begin the induction step we make an initial choice of points $p_{ij}$ in the components of $S_t \cap S[a_i,b_i]$ disjoint from the base.  To achieve disjointness from the previously chosen $p_{kl}$ we will use a simple ``scattering" trick to replace each $p_{ij}$ by a finite number of nearby points, using the following easy fact whose proof we leave as an exercise:

\proclaim Lemma 4.2. Let $K$ be a closed set in $D^p\times D^q$ intersecting each slice $D^p\times\{t\}$ in a finite set, with $p>0$. Then there exist finitely many distinct points $p_r\in D^p$ with corresponding open disks $V_r$ covering $D^q$ such that $K$ is disjoint from $\{p_r\}\times \overline V_r$ for each $r$.

\noindent
We apply this once for each point $p_{ij}$, with $D^p$ a neighborhood of $p_{ij}$ in $S_t$, $D^q=\overline U_i$, and $K$ the union of the previously chosen points $p_{kl}$ for $t$-values in $\overline U_i$.  The result is that each $p_{ij}$ is replaced by a number of points $p_{ijr}$ with the desired disjointness.  For convenience we relabel these points just as $p_{ij}$, each $p_{ij}$ living over a subball $\overline U_{ij}$ of $\overline U_i$.  This gives the induction step for choosing the points $p_{ij}$.

Having chosen all the pairs $p_{ij},q_{ij}$ to be disjoint for each $t$, we can deform the family $S_t$ to create tubes joining these pairs, damping the deformation for the tube from $p_{ij}$ to $q_{ij}$ down to zero near the boundary of $\overline U_{ij}$, with this damping being done in such a way that the full isotopies take place in slightly shrunk versions of $U_{ij}$ that still cover $D^q$.  
Since the ambient space is $\R^\infty$, we can choose all the tubes to be disjoint at all times.  Over $\bdy D^q$ we can choose the tubes to lie outside the slices $S(a_i)$ and the slabs between these slices since the points of $S_t$ in these slices and slabs already lie in the base component of $S_t$, and we arranged in step (1) that these slices and slabs lie inside $S[-1,1]$, so we can choose the tubes to go to either $+\infty$ or $-\infty$ so as to miss these slices and slabs.  This guarantees that the deformation of $f$ that we have constructed induces a corresponding deformation of $\tilde f$ over $\bdy D^q$. 

\smallskip\noindent
(3) \enskip The final step is to make each $S_t$ intersect at least one slice $S(a)$ in $S[-1,1]$ transversely in a connected set, a single circle.  The idea will be to insert tubes that lie in a neighborhood of $S(a)$ so as to join different circles of $S_t\cap S(a)$ together, as shown in the figure below.  Doing this for parametrized families will require a more refined variant of the procedure we used in the preceding step to connect different components of $S_t$ together.
\vskip14pt
\centerline{\epsfbox{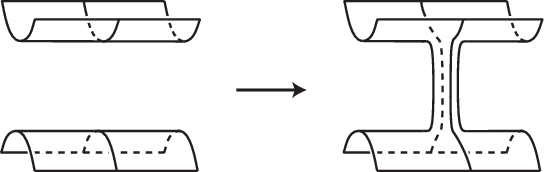}}
\vskip8pt
\noindent
Suppose $S_t$ meets a slice $S(a)$ in $S[-1,1]$ transversely in more than one circle.  For each circle $C_j$ of $S_t\cap S(a)$ that does not intersect the base choose a path $\alpha_j$ in $S_t -B$ from a point $p_j\in C_j$ to a point $q_j$ near $B$.  Step (2) guarantees that such a path exists.  We can assume that $\alpha_j$ is an embedded arc.  (It will not matter if $\alpha_j$ intersects some $C_k$'s at interior points of $\alpha_j$.)  We extend the end of $\alpha_j$ near $B$ by an arc to $-\infty$ parallel to $B$, still calling the extended arc $\alpha_j$. 
We also insert corresponding disjoint arcs $\beta_j$ parallel to $B$ and closer to it, extending to $+\infty$ as well as $-\infty$. 

If the arcs $\alpha_j$ are in general position, they will intersect only in interior points that lie in the sections not parallel to the baseline, and then we can eliminate any intersections of $\alpha_j$ with other $\alpha_k$'s by pushing these $\alpha_k$'s off the end of $\alpha_j$ at $p_j$ by an isotopy $S_t\to S_t$ supported near $\alpha_j$.  This does not introduce any intersections among the other $\alpha_k$'s since the isotopy of the $\alpha_k$'s is the restriction of an isotopy of $S_t$.  Doing this for each $\alpha_j$ in turn, we can thus take all the $\alpha_j$'s (and $\beta_j$'s) to be disjoint.

We can modify $S_t$ by a surgery operation that attaches a thin tube $T_j$ joining a point $x_j$ of $\alpha_j$ to the point $y_j$ of $\beta_j$ having the same first coordinate as $x_j$.  We embed $T_j$ in $\R^\infty$ so that it lies in a small neighborhood of the slice through $x_j$ and $y_j$.  By letting $x_j$ move along $\alpha_j$ from the end at $-\infty$ to the end at $C_j$ while $y_j$ moves simultaneously along $\beta_j$ we obtain a deformation of $S_t$. Doing this for all the circles $C_j$ at once, we deform $S_t$ to a surface $S'_t$ such that $S'_t\cap S(a)$ is connected.  (There is no problem with making all the tubes $T_j$ disjoint in $\R^\infty$.)  By an isotopy of $S'_t$ supported near $S(a)$ and fixed on $B$ we can then make $S'_t$ coincide with the cylinder $Z$ near $S(a)$. 

For nearby $t$-values the curves $C_j$ in $S_t\cap S(a)$ vary by small isotopies and we can take the arcs $\alpha_j$ and $\beta_j$ to vary by small isotopies as well, staying disjoint for each $t$.  By compactness of $D^q$ we can then choose a cover of $D^q$ by open balls $U_i$ such that over $\overline U_i$ there is a slice $S(a_i)$ in $S[-1,1]$ transverse to $S_t$, with arcs $\alpha_{ij}$ going from $-\infty$ to the components $C_{ij}$ of $S_t\cap S(a_i)$ disjoint from the baseline, along with corresponding arcs $\beta_{ij}$ near the baseline.  For fixed $i$ all the arcs $\alpha_{ij}$ and $\beta_{ij}$ are disjoint for each $t\in \overline U_i$, and our next task will be to make them disjoint as $i$ varies as well.

The arcs $\beta_{ij}$ are parallel to $B$ and can easily be chosen to be disjoint by induction on $i$.  For the arcs $\alpha_{ij}$ we also proceed by induction on $i$, so we assume inductively that we have already made the arcs $\alpha_{kl}$ with $k<i$ disjoint for each $t$.  We will show in the next paragraph how we can do a preliminary adjustment so that for each $t\in \overline U_i$, no arc $\alpha_{ij}$ passes through the endpoint $p_{kl}$ of another $\alpha_{kl}$ with $k<i$. Assuming this has been done, any intersections of $\alpha_{ij}$ with arcs $\alpha_{kl}$ for $k<i$ will involve only interior points of $\alpha_{kl}$.  These intersections can be eliminated by the same sort of procedure as before, pushing the arcs $\alpha_{kl}$ off the end of $\alpha_{ij}$ at $p_{ij}$ by an isotopy $S_t \to S_t$ supported near $\alpha_{ij}$, damping this isotopy down to zero outside $U_i$.  

\parshape=11  0pt\hsize  0pt\hsize  0pt\hsize  0pt\hsize  0pt\hsize  0pt\hsize 0pt\hsize 0pt\hsize 0pt\hsize 0pt\hsize 0pt.63\hsize 
\ignorespaces
The preliminary adjustment so that arcs $\alpha_{ij}$ are disjoint from the endpoints $p_{kl}$ of arcs $\alpha_{kl}$ with $k<i$ can be done by a variant of the scattering trick used in step (2) above.  First thicken the arcs $\alpha_{ij}$ to narrow bands of arcs parallel to $\alpha_{ij}$.  Then for a fixed $t$ all but finitely many of these parallel arcs will be disjoint from the finitely many endpoints $p_{kl}$ with $k<i$.  Choose one of these arcs disjoint from $p_{kl}$'s for each $j$.  These choices extend continuously for nearby values of $t$ as well.  Thus we can cover $\overline U_i$ by finitely many neighborhoods $U_{ir}$ with new choices of arcs $\alpha_{ijr}$ in each of these neighborhoods, such that $\alpha_{ijr}$ satisfies the disjointness condition we are trying to achieve. 
We can choose the arcs $\alpha_{ijr}$ to be disjoint from each other for each $t\in \overline U_i$ since we are free to choose $\alpha_{ijr}$ from an open set of the arcs parallel to $\alpha_{ij}$. We choose corresponding slices $S(a_{ir})$ near $S(a_i)$ and connect $\alpha_{ijr}$ to the corresponding circle $C_{ijr}$ of $S_t\cap S(a_{ir})$ as in the figure. 
\vadjust{\hfill\smash{\lower 44pt\llap{\epsfbox{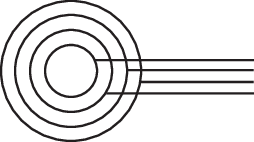}}}}
 To complete the induction step we then replace $U_i$, $S(a_i)$, $C_{ij}$, and $\alpha_{ij}$ by the collections $U_{ir}$, $S(a_{ir})$, $C_{ijr}$, and $\alpha_{ijr}$ and then relabel to eliminate the extra subscripts $r$.

Having all the arcs $\alpha_{ij}$ and $\beta_{ij}$ disjoint, we can use these to construct a well-defined deformation of the family $S_t$ to create tubes $T_{ij}$ over $U_i$ making $S_t$ agree with the cylinder near $S(a_i)$ as described earlier, damping this deformation down near the boundary of $U_i$ as usual.  For the resulting family $S'_t$ there is one problem remaining, however.  Because the deformations for each $U_i$ are damped off near the boundary of $U_i$, the tubes $T_{ij}$ may intersect the slices $S(a_k)$ for other $U_k$'s, perhaps destroying the connectedness of the intersections $S'_t\cap S(a_k)$.  To avoid this problem we first make sure the tubes $T_{ij}$ are very thin while they move to their final destination.  Then there will always be plenty of slices $S(a_{kl})$ near $S(a_k)$ that are disjoint from the moving thin tubes, and we use these slices $S(a_{kl})$ instead of the original slices $S(a_k)$.  (This is another instance of the scattering idea.)

Once we have deformed $S_t$ to meet certain slices $S(a_k)$ in single circles, the parts of $S_t$ in the slabs between these slices will be connected since inserting the tubes $T_{ij}$ does not destroy the property that these parts of $S_t$ lie in the base component of $S_t$, from step (2).  Thus we have constructed a homotopy of the family $f(t)\in \Cinfone_b$ to a family $g(t)$ having a lift $\tilde g(t)\in \Cinfone_s$ obtained by choosing weights for the slices $S(a_k)$.  Over $\bdy D^q$ the slices $S(a_i)$ of $\tilde f(t)$ all lie in the slab $S[-1,1]$ by step (1), so in step (3) there is no need to modify $S_t$ for $t\in\bdy D^q$. \qed 

The argument in the preceding proof can be applied to manifolds of higher dimension $d>2$, but it proves a weaker result, namely that $\Cinfone$ is homotopy equivalent not to the classifying space of a monoid, but to the classifying space of a topological category, a category whose objects are smooth closed connected oriented manifolds of dimension $d-1$ embedded in $\R^\infty$ and whose morphisms are connected oriented cobordisms of dimension $d$ embedded in slabs $[0,a] \times\R^\infty$.  
\medskip

\page

\section 5. Some Variants

\subsect Nonorientable Surfaces

The restriction to oriented surfaces can easily be dropped.  The statement is then that for the standard nonorientable surface $N_\infty$ of infinite genus there is an isomorphism 
$$
H_*(B\Diff_c(N_\infty))\iso H_*(\Omega^\infty_0\overline{AG}^+_{\infty,2})
$$
where $\overline{AG}$ is the version of $AG$ without chosen orientations on the affine planes.  Thus $\overline{AG}^+_{n,2}$ is the Thom space of a bundle over the Grassmann manifold of unoriented $2$-planes in $\R^n$, namely the bundle of vectors normal to these $2$-planes.  In the nonorientable case the monoid $\pi_0\M^\infty$ is not $\Z_{\ge 0}$ but something slightly more complicated, corresponding to diffeomorphism classes of closed connected surfaces under connected sum.  The group completion of this monoid is still $\Z$, however.  The theorem can be restated in terms of mapping class groups of nonorientable surfaces since these satisfy homology stability by [W1] and the Earle-Eells theorem applies also to nonorientable surfaces.

\subsect Punctured Surfaces

Surfaces with punctures can be treated in the same way, provided that one stabilizes with respect to both genus and number of punctures.  Viewing the punctures as distinguished points on a surface rather than deleted points, let $\Diff_c(S_\infty,P)$ be the group of compactly supported diffeomorphisms of the infinite genus surface $S_\infty$ that leave an infinite discrete closed set $P$ in $S_\infty$ invariant, perhaps permuting finitely many of the points of $P$.  The space $AG_{n,2}$ is enlarged to a space $A^*G_{n,2}$ of oriented affine $2$-planes in $\R^n$ with at most one distinguished point in the $2$-plane, where the case of no distinguished point is regarded as the limiting case of a distinguished point that approaches infinity in the $2$-plane.  The statement of the theorem is that there is an isomorphism
$$
H_*(B\Diff_c(S_\infty,P)) \iso H_*(\Omega_0^\infty A^*G^+_{\infty,2})
$$
The space $AG_{n,2}$ is a retract of $A^*G_{n,2}$ by the map that forgets the distinguished point, and this retraction extends to a retraction of one-point compactifications.  The quotient space $A^*G^+_{n,2}/AG^+_{n,2}$ can be identified with the Thom space of the trivial $n$-dimensional vector bundle over $G_{n,2}$ by regarding this bundle as the sum of the two canonical bundles over $G_{n,2}$ consisting of vectors in a given $2$-plane and vectors orthogonal to it;  the vector in the plane gives a distinguished point in the plane, and the vector orthogonal to it tells where to translate the plane.  In the stable homotopy category retractions give wedge sum splittings, so we have an equivalence
$$
\Omega_0^\infty A^*G^+_{\infty,2} \simeq \Omega_0^\infty AG^+_{\infty,2} \times \Omega_0^\infty S^\infty ((G_{\infty,2})_+)
$$
using the fact that the Thom space of a trivial $n$-dimensional bundle over $X$ is the $n$-fold suspension $S^n(X_+)$, the subscript $+$ denoting union with a disjoint basepoint.

This theorem gives information about homology of mapping class groups of finite surfaces since $\Diff_c(S_\infty,P)$ has contractible components and homology stability is known to hold for mapping class groups not just for stabilization with respect to genus but also  with respect to the number of punctures [BT], [Han], [HW].

More refined results that cover stabilization with respect to genus for a fixed number of punctures can be found in [BT].

\subsect
The Zero-Dimensional Case

The parts of the proof of the Madsen-Weiss theorem in Sections 2 and 3 that apply to manifolds of arbitrary dimension become somewhat simpler in dimension $0$ and, with only a small extension, suffice to prove the classical Barratt-Priddy-Quillen theorem.  Let us describe the steps and the simplifications.

\smallskip
\noindent
(1)\enskip The space $\C^n$ consists of configurations of discrete sets of points in $\R^n$, and Proposition 2.1 becomes the statement that $\C^n$ is homotopy equivalent to $S^n$, the subspace consisting of configurations with at most one point.  A continuously-varying tubular neighborhood of such a configuration can be obtained by taking balls centered at the points in the configuration, of radius equal to the minimum of $1$ and one-third of the minimum distance from the point to other points in the configuration.  The procedure in the proof of Proposition~2.1 then gives a deformation retraction of $\C_n$ onto $S^n$.

\smallskip
\noindent
(2)\enskip Lemma 3.2 becomes the statement that $\pi_0\Cnk = 0$ for $k>0$, with a trivial proof since a configuration in $\Cnk$ can be pushed to infinity by radial expansion in the $\R^k$ factor from any point not in the projection of the configuration to $\R^k$.

\smallskip
\noindent
(3)\enskip Proposition 3.1 becomes the statement that $B\Mnk \simeq \Cnkplusone$ for $k>0$.  The proof only uses the easier case that $d<k+1$, when the manifold $M_x$ is empty.

\smallskip
\noindent
(4)\enskip The argument for the cases $k>0$ in the proof of Proposition 3.1 works just as well for $k=0$ to give an equivalence $B\Mnzero\simeq \Cnone$.  There is no need to replace $\Mnzero$ by a monoid $\M^n$, nor is there a need to let $n$ go to infinity.  Thus the more intricate arguments in Section~4 are completely unnecessary.

\smallskip
\noindent
(5)\enskip The group completion theorem yields an isomorphism $H_*(\cup_g\Mnzero_g)\iso H_*(\Omega_0 B\Mnzero)$ where $\Mnzero_g$ is the component of $\Mnzero$ consisting of configurations with $g$ points. This holds for any $n>0$ including $n=\infty$.

\smallskip
\noindent
(6)\enskip Putting together the previous steps, we get  $H_*(\cup_g\Mnzero_g)\iso H_*(\Omega_0^n S^n)$ for $1\le n \le\infty$.

\smallskip
The case $n=1$ is not interesting since $\M^{1,0}_g$ is obviously contractible, as is $\Omega_0 S^1$.  For $n=2$ the space $\M^{2,0}_g$ is homotopy equivalent to the space of configurations of $g$ unordered points in the plane, a classifying space for the braid group $B_g$ for braids with $g$ strands.  The other case when $\Mnzero_g$ is a $K(\pi,1)$ is when $n=\infty$ and $\M^{\infty,0}_g$ is a classifying space $B\Sigma_g$ for the symmetric group $\Sigma_g$.  Thus we have proved the Barratt-Priddy-Quillen theorem and a theorem of Fred Cohen:

\proclaim Theorem.  $H_*(\cup_g\Sigma_g)\iso H_*(\Omega_0^\infty S^\infty)$ and $H_*(\cup_g B_g)\iso H_*(\Omega_0^2 S^2)$.

As an easy extension of these results, we can consider configurations in $\Cn$ with labels in a fixed space $X$ attached to each of the points in a configuration.  Such labeled configurations form a space $\Cn(X)$ topologized to allow the labels to vary continuously over $X$.  There is no difficulty in extending the preceding arguments to the context of labeled configurations since the only deformations of configurations that were used involved pushing points to infinity, and this can be done for labeled points just by keeping the labels unchanged.  The space $\Cn(X)$ has the homotopy type of the subspace of single-point or empty labeled configurations, so this is $(S^n\times X)/(\infty\times X)$ which is the same as $S^n(X_+)$, the $n$-fold reduced suspension of $X$ with a disjoint basepoint adjoined.  Assuming for simplicity that $\pi_0X=0$ so that the path-components of $\Mnzero(X)$ are the subspaces $\Mnzero_g(X)$ with $g$ labeled points, we conclude that there are isomorphisms
$$
\hbox{$H_*(\cup_g\Mnzero_g(X))\iso H_*(\Omega_0^n S^n(X_+))$ for $1\le n \le\infty$}
$$
Ignoring labelings gives a fiber bundle $X^g\To \Mnzero_g(X)\To \Mnzero_g$.  For example if we take $n=\infty$ and $X=BG$ for a discrete group $G$ then $\M^{\infty,0}_g(BG)$ is a classifying space for the wreath product $G\wr \Sigma_g$ and we obtain isomorphisms
$$
H_*(\cup_g (G\wr \Sigma_g))\iso H_*(\Omega^\infty_0 S^\infty(BG_+))
$$
The homology of $G\wr \Sigma_g$ is known to stabilize with $g$, as shown in [HW] for example.  One could also take $n=2$, replacing symmetric groups by braid groups, to obtain isomorphisms
$$
H_*(\cup_g (G\wr B_g))\iso H_*(\Omega^2_0 S^2(BG_+))
$$
where the homology on the left again stabilizes by [HW].

\page

\vskip-15pt

\section Appendix A. Classifying Spaces for Diffeomorphism Groups

The most classical of classifying spaces are the classifying spaces for vector bundles.  The classifying space for real vector bundles of dimension $k$ is the Grassmannian $G_{\infty,k}$ of $k$-dimensional vector subspaces of $\R^\infty$, the direct limit of the Grassmann manifolds $G_{n,k}$ of $k$-dimensional vector subspaces of $\R^n$.  (This notation is not quite consistent with the notation earlier in the paper, where $G_{\infty,k}$ meant the $2$-sheeted cover consisting of oriented $k$-dimensional vector subspaces of $\R^\infty$.)  There is a canonical vector bundle $E_{\infty,k}$ over $G_{\infty,k}$ consisting of the pairs $(v,P)\in\R^\infty\times G_{\infty,k}$ with $v\in P$.  This is a universal $k$-dimensional vector bundle in the sense that every vector bundle $E\to X$ is induced from the universal bundle by some map $X\to G_{\infty,k}$ which is unique up to homotopy.  Some mild restrictions on $X$ are needed for this, for example that $X$ is paracompact.  

If we shift our point of view from vector spaces to automorphisms of vector spaces, there is an associated fiber bundle over $G_{\infty,k}$ whose fibers are $GL(k,\R)$ rather than $\R^k$. This bundle is often written as $EGL(k,\R)\to BGL(k,\R)$ where $BGL(k,\R)$ is just another notation for $G_{\infty,k}$ and $EGL(k,\R)$ is the space of linear embeddings $f\:\R^k\to \R^\infty$, with the projection $EGL(k.\R)\to BGL(k,\R)$ sending an embedding $f$ to its image $f(\R^k)$.  Since a linear embedding $f$ is determined by where it sends the standard basis of $\R^k$, one could also describe $EGL(k,\R)$ as the Stiefel manifold of $k$-tuples of linearly independent vectors in $\R^\infty$.  The condition of linear independence can be strengthened to orthonormality, which amounts to requiring the embeddings $f$ to be isometric embeddings, and then $GL(k,\R)$ is replaced by the orthogonal group $O(k)$ and one has a fiber bundle $EO(k)\to BO(k)$ with fiber $O(k)$.  This does not affect the homotopy types of the spaces, and in fact $BO(k)$ is the same space as $BGL(k,\R)$, namely $G_{\infty,k}$.  

A key feature of $EGL(k,\R)$ and $EO(k)$ is that they are contractible.  It is an elementary fact (see [H1], Proposition~4.66) that whenever one has a fiber bundle or fibration $F \to E \to B$ with $E$ contractible, then there is a (weak) homotopy equivalence $F \to \Omega B$.  Thus $O(k)\simeq \Omega BO(k) =\Omega G_{\infty,k}$.

Entirely analogous considerations hold for smooth fiber bundles with fiber a smooth compact manifold $M$.  A classifying space for such bundles is the space $B\Diff(M)=\C(M,\R^\infty)$ of smooth submanifolds of $\R^\infty$ diffeomorphic to $M$. (The notation $\C(M,\R^\infty)$ conflicts with that used earlier in the paper since we are ignoring orientations now.)  The topology on $\C(M,\R^\infty)$ is the direct limit of its subspaces $\C(M,\R^n)$ which are in turn given the quotient topology obtained by regarding $\C(M,\R^n)$ as the orbit space of the space of smooth embeddings $M\to \R^n$ (with the usual $C^\infty$ topology) under the action of $\Diff(M)$ by composition.  

There is a canonical bundle $\E(M,\R^\infty)\to\C(M,\R^\infty)$ with fiber $M$, where $\E(M,\R^\infty)$ consists of pairs $(v,P)$ in $\R^\infty \times \C(M,\R^\infty)$ with $v \in P$.  To see that this bundle is locally a product, consider a given $P\subset \R^n$ diffeomorphic to $M$ with an open tubular neighborhood $N$ that is identified with the normal bundle to $P$ in $\R^n$.  Then the sections of this bundle form a neighborhood of $P$ in $\C(M,\R^n)$, and projecting these sections onto the $0$-section $P$ gives a local product structure for the projection $\E(M,\R^n)\to\C(M,\R^n)$, compatibly for increasing $n$, hence also for the direct limit $\E(M,\R^\infty)\to\C(M,\R^\infty)$.

The bundle $\E(M,\R^\infty)\to\C(M,\R^\infty)$ is universal for smooth bundles with fiber $M$ over paracompact base spaces, by essentially the same argument that shows the vector bundle $E_{\infty,k}\to G_{\infty,k}$ is universal; see for example the first chapter of [H2].  Realizing a given bundle $M\to E \to B$ as a pullback of the universal bundle is equivalent to finding a map $E\to\R^\infty$ which is a smooth embedding on each fiber.  Locally in the base $B$ such maps exist by combining a local projection onto the fiber $M$ with a fixed embedding $M\to\R^n$.  These local maps to $\R^n$ that are embeddings on fibers can then be combined to a global map $E\to\R^\infty$ that is an embedding on fibers via a partition of unity argument just as in the vector bundle case.  Thus every bundle $M\to E\to B$ is a pullback of the universal bundle via some map $B\to \C(M,\R^\infty)$.  Uniqueness of this map up to homotopy follows from the uniqueness of the map $E\to\R^\infty$ up to homotopy, which is shown just as for vector bundles by composing with linear embeddings $\R^\infty \to \R^\infty$ onto the odd or even coordinates, using straight-line homotopies.

There is also an associated bundle $\Diff(M)\to E\Diff(M) \to B\Diff(M)$ whose total space $E\Diff(M)$ is the space of smooth embeddings $M\to\R^\infty$. The fiber bundle property can be proved using tubular neighborhoods and sections as before.  The space $E\Diff(M)$ is contractible by pushing each embedding $M\to\R^\infty$ into the odd coordinates by composing with a linear isotopy of $\R^\infty$ into the odd coordinates, then taking a linear isotopy to a fixed embedding of $M$ into the even coordinates. 

Since we have a fiber bundle $\Diff(M)\to E\Diff(M) \to B\Diff(M)$ with contractible total space, it follows that $\Diff(M)$ is weakly equivalent to $\Omega B\Diff(M)$.  In particular, $B\Diff(M)$ is a $K(\pi,1)$ for the mapping class group of $M$ if the components of $\Diff(M)$ are contractible.

\page

\section Appendix B. The Earle--Eells Theorem

Let $S$ be a compact connected surface, not necessarily orientable and possibly with boundary.  In this appendix we use the notation $\Diff(S)$ for the group of diffeomorphisms of $S$ that restrict to the identity on $\bdy S$, with the $C^\infty$ topology as usual.

\proclaim Theorem B.1. The components of $\Diff(S)$ are contractible except when $S$ is the sphere, projective plane, torus, or Klein bottle. 

In the exceptional cases the components of $\Diff(S)$ are homotopy equivalent to $SO(3)$ for the sphere and projective plane, $S^1\times S^1$ for the torus, and $S^1$ for the Klein bottle.  These statements could be proved by mild extensions of the arguments below.

The theorem was first proved in [EE] and [ES] for orientable surfaces. We will give an exposition of Gramain's proof in [Gr].  This is based on a proof by Cerf for the case $S=D^2$ in [C1], later republished in the appendix of the often-cited volume [C2].  This special case already contains the key trick that makes the proof so simple.  In the write-up below this trick occurs at the very end.

\medskip
\pf Proof. There will be three main steps.

\smallskip\noindent
(I)\enskip {\it Reduction to the case of nonempty boundary.}  

\noindent
Assume $\bdy S$ is empty.  By evaluating diffeomorphisms $S\to S$ at a basepoint $x_0\in S$ we obtain a fibration
$$
\Diff(S,x_0) \To \Diff(S)\To S
$$
whose fiber $\Diff(S,x_0)$ consists of diffeomorphisms fixing $x_0$.  From the long exact sequence of homotopy groups for this fibration we obtain isomorphisms $\pi_i\Diff(S)\iso\pi_i\Diff(S,x_0)$ for $i>1$ since $\pi_i(S)=0$ for $i>1$ from the assumption that $S\ne S^2,P^2$.  We can also deduce that $\pi_1\Diff(S)\iso\pi_1\Diff(S,x_0)$ by looking at the end of the long exact sequence
$$
0\To\pi_1\Diff(S,x_0)\To\pi_1\Diff(S)\To \pi_1(S,x_0) \mapright \bdy \pi_0\Diff(S,x_0)
$$
where it suffices to show that the final map $\bdy$ is injective.  To verify injectivity, recall from the construction of the long exact sequence that $\bdy$ sends the homotopy class $[\gamma]$ of a loop $\gamma$ in $S$ at the basepoint $x_0$ to the isotopy class $[f]$ of the diffeomorphism $f\:(S,x_0)\to(S,x_0)$ obtained by dragging $x_0$ around $\gamma$ and extending this to an isotopy of $S$ ending at $f$.  The automorphism of $\pi_1(S,x_0)$ induced by $f$ is conjugation by $\gamma$, so the composition $\pi_1(S,x_0)\mapright \bdy \pi_0\Diff(S,x_0)\To \Aut(\pi_1(S,x_0))$ is injective since $\pi_1(S,x_0)$ has trivial center for $S$ as in the theorem.  Hence $\bdy$ is injective. 

Next we look at another fibration, this one obtained by evaluating diffeomorphisms of $(S,x_0)$ on a closed disk $D$ which is a neighborhood of $x_0$.
$$
\Diff(S,D)\To\Diff(S,x_0)\To \Emb((D,x_0),(S,x_0))
$$
The fiber $\Diff(S,D)$ consists of diffeomorphisms fixing $D$ and the base is the space of smooth embeddings $D\to S$ fixing $x_0$.  We have $\Diff(S,D)\simeq\Diff(S_0)$ for $S_0=S-\int(D)$.  There is a projection $\Emb((D,x_0),(S,x_0))\to GL(2,\R)$ obtained by taking the derivative at $x_0$, and this is a homotopy equivalence as is true for manifolds of any dimension by a very special case of tubular neighborhood theory.  Since $GL(2,\R)\simeq O(2)$ we deduce that $\pi_i\Diff(S,x_0)\iso\pi_i\Diff(S,D)$ for $i>1$, and we can again extend this to the case $i=1$ by examining the exact sequence:
$$
0\To\pi_1\Diff(S,D)\To\pi_1\Diff(S,x_0)\To \pi_1\Emb((D,x_0),(S,x_0)) \mapright \bdy \pi_0\Diff(S,D)
$$
Here we have $\pi_1\Emb((D,x_0),(S,x_0))\iso\Z$ generated by a full rotation of $D$ about $x_0$.  The map $\bdy$ sends this rotation to the diffeomorphism obtained by extending this rotation to a diffeomorphism of $S$, which we can take to be a Dehn twist about a curve parallel to $\bdy D$.  This twist and all its nonzero powers induce nontrivial inner automorphisms of the free group $\pi_1(S_0,y_0)$ for a basepoint $y_0\in\bdy S_0$, so injectivity of $\bdy$ follows.

Thus we have shown that $\pi_i\Diff(S)\iso\pi_i\Diff(S,x_0)\iso\pi_i\Diff(S_0)$ for all $i>0$, which gives the reduction to the case that $\bdy S$ is nonempty.

\medskip\noindent
(II)\enskip {\it Reduction to contractibility of certain spaces of arcs.} 

\noindent
We assume that $\bdy S$ is nonempty from now on.  Consider smooth embeddings $I\to S$ joining two given points $p$ and $q$ in $\bdy S$ and whose images are proper arcs in $S$, intersecting $\bdy S$ only in their endpoints.  Fixing one such embedding $\alpha$, let $A(S,\alpha)$ be the space of all such embeddings $I\to S$ isotopic to $\alpha$ fixing the endpoints.  In step (III) we will show that $A(S,\alpha)$ is contractible, and let us now show that this implies that the path-components of $\Diff(S)$ are contractible.  

Since all the path-components are homeomorphic, it suffices to show this for the path-component $\Diff_0(S)$ of the identity map.  

Evaluation of diffeomorphisms on $\alpha(I)$ gives a fibration
$$
\Diff_0(S,\alpha)\To \Diff_0(S)\To A(S,\alpha)
$$
whose fiber consists of diffeomorphisms in $\Diff_0(S)$ restricting to the identity on $\alpha(I)$.  Under the assumption that $A(S,\alpha)$ is contractible, we have $\pi_0\Diff_0(S,\alpha)\iso\pi_0\Diff_0(S)$ from the long exact sequence of homotopy groups, so $\Diff_0(S,\alpha)$ is the identity component of $\Diff(S,\alpha)$, the diffeomorphisms of $S$ fixing $\bdy S$ and $\alpha(I)$.  Thus $\Diff_0(S,\alpha)\simeq\Diff_0(S')$ where $S'$ is obtained from $S$ by cutting along $\alpha$.  From the long exact sequence we conclude that $\Diff_0(S)\simeq \Diff_0(S')$.  

We can reduce $S$ to a disk by a finite sequence of cuts along nonseparating arcs, so by induction on $-\chi(S)$ the problem is reduced to showing that $\Diff_0(D^2)$ is contractible.  To show this we consider the fibration
$$
\Diff(D^2_+)\To\Emb(D^2_+,D^2)\To A(D^2,\alpha)
$$
where $D^2_+$ is the upper half of $D^2$ and $\Emb(D^2_+,D^2)$ is the space of embeddings $D^2_+\to D^2$ fixing $D^2_+\cap \bdy D^2$ and taking the rest of $\bdy D^2_+$ to the interior of $D^2$.  Restriction to $D^1\subset \bdy D^2_+$ gives the map to $A(D^2,\alpha)$ with $\alpha$ the inclusion $D^1\to D^2$.  The space $\Emb(D^2_+,D^2)$ is contractible, by the standard argument that shows the space of embeddings of a collar on the boundary of a manifold is contractible.  So the map to $A(D^2,\alpha)$ does indeed produce embeddings isotopic to $\alpha$.  From the long exact sequence we see that contractibility of $A(D^2,\alpha)$ implies contractibility of $\Diff(D^2_+)$ and hence of $\Diff(D^2)$. (In particular this shows that $\Diff_0(D^2)=\Diff(D^2)$.)

\medskip\noindent
(III)\enskip {\it Contractibility of the arc spaces.}  

\noindent
It remains to show the spaces $A(S,\alpha)$ are contractible.  The easier case is that $\alpha$ connects two different components of $\bdy S$ so we do this case first.  Let $T$ be the surface obtained from $S$ by filling in the component of $\bdy S$ at one end of $\alpha$ with a disk, say the end at $q$.  Consider the fibration
$$
\Emb(I,S)\To\Emb(I\cup D^2,T)\To\Emb(D^2,T-\bdy T)
$$
where $I\cup D^2$ is formed by attaching $1\in I$ to a point in $\bdy D^2$, and $\Emb(I\cup D^2,T)$ is the space of embeddings $f\:I\cup D^2\to T$ taking $0$ to $p$ and the rest of $I\cup D^2$ to the interior of $T$.  We assume also that $f$ preserves orientation on $D^2$, where we use $f(I)$ to transport an orientation of $T$ at $p$ to an orientation at $f(1)$.  It is not hard to see that $\Emb(I\cup D^2,T)$ is contractible.  For example, by restricting embeddings $f$ to $I$ we get a fibration over the space of embeddings $I\to T$ taking $0$ to $p$ and the rest of $I$ to the interior of $T$, obviously a contractible space, and the fibers are also evidently contractible.

The base space $\Emb(D^2,T-\bdy T)$ in the preceding displayed fibration has $\pi_i=0$ for $i>0$.  (One can see this by fibering over $T-\bdy T$ by evaluating embeddings $D^2\to T-\bdy T$ at $0$.)  From the long exact sequence of homotopy groups we conclude that the fiber $\Emb(I,S)$ has contractible components.  One of these components is $A(S,\alpha)$, so this takes care of the case that $p$ and $q$ lie in different components of $\bdy S$. 

\smallskip
The harder case is when $\alpha$ joins a component of $\bdy S$ to itself.  (This is the one place in the proof where some cleverness is required, and where the fact that we are dealing with surfaces plays a crucial role.)  In this case let $T$ be obtained from $S$ by attaching a rectangle $R$ by identifying two opposite edges of $\bdy R$ with arcs in $\bdy S$ close to $p$ on either side of $p$. (We can view $R$ as a $1$-handle.)  The effect of attaching $R$ is that $p$ and $q$ now lie in distinct components of $\bdy T$.  Let $\beta$ be an arc in $R$ joining the two nonattaching edges.  Cutting $T$ along $\beta$ produces a surface $U$ homeomorphic to $S$.  

\vskip6pt
\centerline {\epsfbox{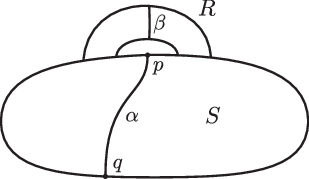}}

\noindent
It will suffice to show that the maps $\pi_iA(T-\beta,\alpha)\to\pi_iA(T,\alpha)$ are injective for all $i$ since the latter groups are trivial by the case already proved.  

We prove injectivity of $\pi_iA(T-\beta,\alpha)\to\pi_iA(T,\alpha)$ by lifting to a certain covering space $\tilde T$ of $T$.  Since $T\simeq S\vee S^1$ we have $\pi_1T \iso \pi_1S *\Z$, and we let $\tilde T$ be the covering space corresponding to the subgroup $\pi_1S$ of $\pi_1S * \Z$.  Explicitly, $\tilde T$ is built from $U$ by attaching two copies of the universal cover $\tilde U$ of $U$ to $U$ along the two copies of $\beta$ in $\bdy U$.  

\vskip10pt
\centerline {\epsfbox{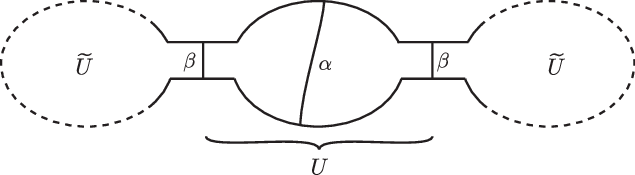}}

All arcs in $A(T,\alpha)$ are isotopic to $\alpha$ so they have unique lifts to $\tilde T$ isotopic to the lifted copy of $\alpha$ in $U \subset \tilde T$.  These lifts form a subspace $\tilde A(T,\alpha)$ of $A(\tilde T,\alpha)$ homeomorphic to $A(T,\alpha)$ and we can identify the inclusion $A(T-\beta,\alpha)\incl A(T,\alpha)$ with the first of the two inclusions
$$
A(U,\alpha)\incl \tilde A(T,\alpha) \incl A(\tilde T,\alpha)
$$
It will suffice to show that the composition $i\: A(U,\alpha)\To  A(\tilde T,\alpha)$ of the two inclusions induces injections on homotopy groups.  We will do this by constructing a map in the opposite direction $r\: A(\tilde T,\alpha)\To A(U,\alpha)$ such that $ri$ is homotopic to the identity.  

The interior of $\tilde U$ is diffeomorphic to $\R^2$ and $\tilde U - \bdy \tilde T$ is diffeomorphic to $\R \times [0,\infty)$.  Let $\tilde T_0$ be $\tilde T$ with the part of $\bdy \tilde T$ lying in the two copies of $\tilde U$ deleted.  Via the two products $\R \times [0,\infty)$ in $\tilde T_0$ we can isotope $\tilde T_0$ into $U$.  The final map in this isotopy induces the map $r$ and the restriction of the isotopy to $U$ gives a homotopy from $ri$ to the identity.  This finishes the proof.   \qed

In the special case $S=D^2$ this argument shows that the space of smooth proper arcs in $D^2$ with fixed endpoints has contractible components, but it does not show this space is path-connected, which is equivalent to the Schoenflies theorem in the smooth category.

\section Appendix C. Rational Homology

The main result we prove here is:

\proclaim Theorem C.1. $H^*(\Omega^\infty_0 AG^+_{\infty,2};\Q)$ is a polynomial algebra $\Q[x_2,x_4,x_6,\cdots\,]$ on an infinite sequence of even-dimensional generators $x_{2i}\in H^{2i}$. 

As we will see, the method of proof is fairly general and applies equally well in many other similar situations.

\medskip
\pf Proof. Recall that $AG^+_{n,2}$ is the Thom space of an $(n-2)$-dimensional vector bundle over the Grassmann manifold $G_{n,2}$ of oriented $2$-planes in $\R^n$ containing the origin.  This is an oriented vector bundle, so there are Thom isomorphisms $\tilde H^i(AG^+_{n,2})\approx H^{i-(n-2)}(G_{n,2})$ for all $i$.  

The cohomology of $G_{n,2}$ is the same as for $\CP^\infty$ in dimensions up to approximately $n$, as one can see from the fiber bundle $S^1\to V_{n,2}\to G_{n,2}$ involving the Stiefel manifold $V_{n,2}$ of orthonormal $2$-frames in $\R^n$.  Since $V_{n,2}$ fibers over $S^{n-1}$ with fiber $S^{n-2}$, it is $(n-3)$-connected, so from the long exact sequence of homotopy groups for the bundle $S^1\to V_{n,2}\to G_{n,2}$ we see that $G_{n,2}$ is a $K(\Z,2)$ up to dimension approximately $n$.  Thus there is a map $G_{n,2}\to\CP^{\infty}$ inducing an isomorphism on homotopy groups up to dimension about $n$, hence also on homology and cohomology.  This map is essentially just the inclusion $G_{n,2}\subset G_{\infty,2}$ since $G_{\infty,2}$ is also a $K(\Z,2)$ via the bundle $S^1\to V_{\infty,2}\to G_{\infty,2}$ with $V_{\infty,2}$ contractible.  Generators for the groups $H^{2i}(G_{\infty,2})$ are the powers $e^i$ of the Euler class of the universal bundle, so the restrictions of these classes $1,e,e^2,\cdots$ to $G_{n,2}$ give an additive basis for $H^*(G_{n,2})$ up to dimension approximately $ n$.  

Via the Thom isomorphism the classes $1,e,e^2,\cdots$  give generators for $\tilde H^*(AG^+_{n,2})$ up to dimension around $2n$.  Since cohomology is represented by maps to Eilenberg-MacLane spaces, these generators for $\tilde H^*(AG^+_{n,2})$ correspond to a map 
$$
f\:AG^+_{n,2}\To K(\Z,n-2)\times K(\Z,n) \times K(\Z,n+2) \times \cdots
$$
Now we switch to rational coefficients to take advantage of the fact that the rational cohomology of Eilenberg-MacLane spaces $K(\Z,k)$ is rather simple, namely a polynomial algebra $\Q[x]$ when $k$ is even and an exterior algebra $\Lambda_{\smQ}[x]$ when $k$ is odd, with $x$ lying in $H^k(K(\Z,k);\Q)$ in both cases.  This standard calculation is done by induction on $k$ using the path fibration $K(\Z,k-1)\to P\to K(\Z,k)$ with contractible total space $P$.  For the actual calculation one can use either the Serre spectral sequence or the Gysin and Wang exact sequences with $\Q$ coefficients.  From the K\"unneth formula it follows that the map $f$ is an isomorphism on rational cohomology in dimensions up to about $2n$ since the first nontrivial product in the cohomology of the space on the right occurs in dimension $2n-4$.  

From the universal coefficient theorem it follows that $f$ is also an isomorphism on rational homology in the same range.  We wish to deduce from this that $f$ is an isomorphism on rational homotopy groups $\pi_*\otimes\Q$ in this range.  The quickest explanation for this fact is via $\Q$-localization.  The map $f$ induces a map $f_{\smQ}$ of the $\Q$-localizations of the two spaces, and $\Q$-localizations are characterized by the fact that their integer and rational homology groups are isomorphic, as are their homotopy groups and rational homotopy groups.  Thus $f_{\smQ}$ induces isomorphisms on integer homology up to dimension roughly $2n$, so by the relative Hurewicz theorem it also induces isomorphisms on homotopy groups in this range, so $f$ induces isomorphisms on rational homotopy groups in this range.  (Note that the spaces involved here are all simply-connected if $n$ is large enough, so $\Q$-localizations exist and there are no subtleties with the relative Hurewicz theorem.)

Now we apply the loopspace functor repeatedly, using the fact that $\Omega(X\times Y)=\Omega X\times\Omega Y$.  After applying it $n-2$ times we have a map $$\Omega^{n-2}AG^+_{n,2}\To K(\Z,0)\times K(\Z,2) \times K(\Z,4) \times \cdots$$ inducing an isomorphism on rational homotopy groups up to dimension roughly $n$, since the range was up to about $2n$ before.  Applying $\Omega$ one more time involves only the basepoint components of both spaces, so the factor $K(\Z,0)$ will disappear, and we will have a map $$\Omega^{n-1}AG^+_{n,2}\To K(\Z,1)\times K(\Z,3) \times K(\Z,5) \times \cdots$$ so that one final looping will give a map $$\Omega^nAG^+_{n,2}\To K(\Z,0)\times K(\Z,2) \times K(\Z,4) \times \cdots$$
Restricting to a single component of each space, we then have a map
$$
\Omega^n_0 AG^+_{n,2}\To K(\Z,2)\times K(\Z,4) \times K(\Z,6) \times \cdots
$$
This induces isomorphisms on rational homotopy groups up to dimension roughly $n$.  To deduce that this induces isomorphisms on rational homology in this range we can again make use of $\Q$-localization.  The space on the right is simply-connected, but the one on the left might not be.  However, it is an H-space so it has $\pi_1$ abelian with trivial action on all higher homotopy groups, which is sufficient for a $\Q$-localization to exist.  After $\Q$-localization the map is an isomorphism on homotopy groups in a range, so relative homotopy groups are trivial in this range, hence also the relative homology groups, so we get the desired conclusion that the map displayed above is an isomorphism on rational homology up to dimension about $n$.  The same is then true for rational cohomology.  Thus by the K\"unneth formula we have shown that $H^*(\Omega^n_0 AG^+_{n,2};\Q)$ is a polynomial algebra $\Q[x_2,x_4,x_6,\cdots\,]$ up through dimension approximately $n$.  

It remains to check that the cohomology of $\Omega^n_0 AG^+_{n,2}$ stabilizes as $n$ increases, that is, the maps $\Omega^n_0 AG^+_{n,2}\to \Omega^{n+1}_0 AG^+_{n+1,2}$ induce isomorphisms on cohomology in a range of dimensions that goes to infinity with $n$.  For this we can use $\Z$ coefficients.  To start, we consider the map $ AG^+_{n,2}\to \Omega AG^+_{n+1,2}$.  This factors as a composition 
$$
AG^+_{n,2}\To \Omega\Sigma AG^+_{n,2}\To \Omega AG^+_{n+1,2}
$$
where $\Sigma$ denotes reduced suspension and the first map is an instance of the canonical map $X\to \Omega\Sigma X$ defined for any space $X$.  The second map is the looping of the map $\Sigma AG^+_{n,2}\to AG^+_{n+1,2}$ adjoint to $ AG^+_{n,2}\to \Omega AG^+_{n+1,2}$.  The Thom space $AG^+_{n,2}$ is roughly $n$-connected so the first map is an isomorphism on homotopy groups through dimension approximately $2n$ by the Freudenthal suspension theorem, hence also on homology and cohomology in this range.  The map $\Sigma AG^+_{n,2}\to AG^+_{n+1,2}$ is an isomorphism on cohomology through dimension about $2n$ via the Thom isomorphisms for the two spaces (the suspension of a Thom space is the Thom space for the sum of the given bundle with a trivial line bundle), using the fact that the pair $(G_{n+1,2},G_{n,2})$ is roughly $n$-connected via the earlier comparison of these spaces with $G_{\infty,2}$.  Since the map $\Sigma AG^+_{n,2}\to AG^+_{n+1,2}$ is an isomorphism on cohomology up to dimension about $2n$, this is true also for homology (both spaces being finite CW complexes) and so also for homotopy groups.  Looping the map $\Sigma AG^+_{n,2}\to AG^+_{n+1,2}$ loses just one degree of connectivity, so the second of the two composed maps displayed above is also an isomorphism on homotopy groups up to dimension around $2n$. 

Thus the composed map  $ AG^+_{n,2}\to \Omega AG^+_{n+1,2}$ is an isomorphism on homotopy groups up to dimension approximately $2n$, so after applying $\Omega^n_0$ to this map we obtain isomorphisms on homotopy groups to dimension roughly $n$, hence also on homology and cohomology groups.     
 \qed

As this proof makes clear, the first two additive generators for $\tilde H^*(AG^+_{n,2};\Q)$ corresponding to $1$ and $e$ in $H^*(G_{n,2})$ eventually disappear after iterated looping since they corresponded to factors $K(\Z,0)$ that were canceled.  Furthermore, we can see how the remaining additive generators corresponding to the higher powers $e^2,e^3,\cdots$ become multiplicative generators of the polynomial ring $\Q[x_2,x_4,\cdots\,]$ since under iterated looping, more and more of the powers and products of these additive generators fall inside the stable range of dimensions.

It should be clear how the preceding calculation generalizes to give a calculation of $H^*(\Omega^\infty_0 AG^+_{\infty,d};\Q)$ for $d>2$.  The input is the well-known calculation of $H^*(G_{\infty,d};\Q)$ as a polynomial ring on even-dimensional generators, the first approximately $d/2$ Pontryagin classes plus the Euler class when $d$ is even.  An additive basis for $H^*(G_{\infty,d};\Q)$ therefore consists of the monomials in these classes.  After restricting to $G_{n,d}$ for large finite $n$, the Thom isomorphism shifts these additive generators up by $n-d$ dimensions.  Then the $n$-fold looping process shifts these additive generators back down by $n$ dimensions, deleting the generators whose dimensions become negative or $0$.  The surviving additive generators correspond to monomials in $H^*(G_{\infty,d};\Q)$ of dimension greater than $d$, and these become multiplicative generators of $H^*(\Omega^\infty_0 AG^+_{\infty,d};\Q)$ which is a polynomial algebra on even-dimensional generators when $d$ is even and an exterior algebra on odd-dimensional generators when $d$ is odd, since the net shift in dimensions is downward by $d$.

\medskip
To finish this appendix let us extend the rational cohomology calculation to the case of nonorientable surfaces where, as we saw in section~5, the space $AG_{\infty,2}$ is replaced by its quotient $\overline{AG}_{\infty,2}$ obtained by ignoring the orientations of affine $2$-planes in $\R^\infty$.  

\proclaim Theorem C.2. $H^*(\Omega^\infty_0 \overline{AG}^+_{\infty,2};\Q)$ is a polynomial algebra $\Q[y_4,y_8,y_{12},\cdots\,]$ on generators $y_{4i}\in H^{4i}$. 

\pf Proof. The main obstacle to carrying over the previous proof directly is that we can no longer apply the Thom isomorphism to compute $\tilde H^*(\overline{AG}^+_{n,2})$ with $\Z$ or even $\Q$ coefficients since we are now dealing with the Thom space of a nonorientable vector bundle.  Fortunately there is an easy way around this problem by using the classical transfer homomorphism for finite-sheeted covering spaces.  The projection $AG^+_{n,3}\to\overline{AG}^+_{n,2}$ is not a covering space at the compactification point, but this is no longer an issue if we consider relative cohomology for the pairs consisting of the unit disk bundles and their boundary sphere bundles, since these relative cohomology groups are isomorphic to the reduced absolute groups for the Thom spaces.  The basic properties of transfer homomorphisms work equally well for relative cohomology, so we can say that $\tilde H^*(\overline{AG}^+_{n,2};\Q)$ is isomorphic to the subgroup of $\tilde H^*(AG^+_{n,2};\Q)$ consisting of elements invariant under the map induced by the involution of $AG^+_{n,2}$ that reverses the orientation of each $2$-plane.  On $G_{n,2}$ this map sends the Euler class $e$ to its negative, hence $e^i$ is sent to $(-1)^i e^i$.  The Thom class is also sent to its negative.  Since the Thom isomorphism is given by cup product with the Thom class, we conclude that the invariant elements of $\tilde H^*(AG^+_{n,2};\Q)$ in the range we are interested in are the ones corresponding to odd powers of $e$.  These occur in dimensions $n,n+4,n+8,\cdots$.  The class in dimension $n$ corresponding to $e$ itself disappears after $n$-fold looping as before, so we are left with the classes corresponding to $e^3,e^5,e^7,\cdots$, which live in dimensions $4,8,12,\cdots$.   \qed

\vskip-12pt

\section Appendix D. The Group Completion Theorem 

The goal of this appendix is to prove a version of the Group Completion Theorem that suffices for many applications such as the Madsen-Weiss theorem.  Some comments on more general versions are given at the end.

As motivation for the plan of attack, recall an elementary fact about fibrations which is Proposition~4.66 of [H1]:  If $F\to E\to B$ is a fibration with $E$ contractible, then there is a weak homotopy equivalence $F\to \Omega B$.  This is proved by noting that a contraction of $E$ gives rise to a map from the given fibration to the standard path fibration $\Omega B \to PB \to B$ and so by the five lemma the map on fibers $F\to \Omega B$ is a weak homotopy equivalence.  One might hope to apply this fact when $B$ is $B\M$ for a given topological monoid $\M$.  We will construct a contractible space $E\M$ with a projection $E\M\to B\M$ whose fibers are copies of $\M$, but the projection will not generally be a fibration.  However, in some cases it is something a little weaker than a fibration that works just as well, a quasifibration.  

\parshape=13  0pt\hsize  0pt\hsize  0pt\hsize  0pt\hsize  0pt\hsize   0pt.53\hsize 0pt.53\hsize 0pt.53\hsize 0pt.53\hsize 0pt.53\hsize 0pt.53\hsize 0pt.53\hsize 0pt\hsize 
The definition of $E\M$ is similar to that of $B\M$, except that instead of starting with a single vertex, one has a vertex for each element of $\M$.  In fact there is a more general construction that will eventually be needed.  Suppose we are given an action of the monoid $\M$ on a space $X$, this being a right action, so each $m\in\M$ gives a map $X\to X$, $x\mapsto x\cdot m$.  
\vadjust{\hfill\smash{\lower 104pt\llap{\epsfbox{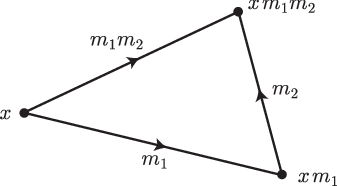}}}}\ignorespaces 
Then we define the space $EX$ to have $X$ as its space of vertices, and with its space of $p$-simplices the product $X \times \M^p$, so there is a $p$-simplex for each $(p+1)$-tuple $(x,m_1,\cdots,m_p)$, with vertices labeled $x,x\cdot m_1, x\cdot m_1m_2,\cdots$ and  edges labeled $m_1,m_2,\cdots$ as in the figure.
Thus $EX$ is a quotient of $\coprod_p \Delta^p\times X\times \M^p$ with identifications similar to those in $B\M$.  In particular, when $X$ is a point then $EX$ is just $B\M$.  For a general $X$, forgetting the vertex labels gives a projection $EX\to B\M$ with fibers $X$.  

\parshape=1   0pt.5\hsize 
A case of special interest is when $X$ is $\M$ acting on itself on the right.  The space $E\M$ is then contractible by letting each $p$-simplex flow to the vertex corresponding to the identity element $e$ in $\M$ by flowing linearly in the $(p+1)$-simplex obtained by coning to the vertex $e$.
\vadjust{\hfill\smash{\lower 0pt\llap{\epsfbox{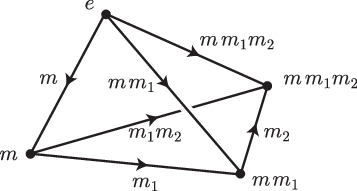}}}}\ignorespaces

\proclaim Lemma D.1.  If $\cdot m\: X\to X$ is a weak homotopy equivalence for each $m\in\M$ then the projection $EX\to B\M$ is a quasifibration.

Recall the definition: a map $p\:E \to B$ is a quasifibration if for each $b\in B$ the map $p_*\:\pi_i(E,p^{-1}(b),x_0)\to\pi_i(B,b)$ is an isomorphism for all $i$ and all $x_0\in p^{-1}(b)$.  These isomorphisms immediately yield a long exact sequence of homotopy groups for a quasifibration that is formally identical to the long exact sequence for a fibration or fiber bundle.

\medskip
\pf Proof of Lemma D.1.  We will use Lemma 4K.3 of [H1] which gives three conditions each of which implies that a map $p\: E\to B$ is a quasifibration.  These conditions can be stated as follows:

\smallskip\noindent
(a)\enskip $B$ is the union of two open sets over each of which $p$ is a quasifibration, and $p$ is also a quasifibration over their intersection.

\smallskip\noindent
(b)\enskip $B$ is the direct limit of an increasing sequence of subspaces $B_n$ over each of which $p$ is a quasifibration.

\smallskip\noindent
(c)\enskip $E$ deforms into a subspace $E'$ via a deformation $F_t\: E\to E$ that covers a deformation of $B$ into a subspace $B'$ such that $p$ restricts to a quasifibration $E'\to B'$ and the map $F_1\: p^{-1}(b)\to p^{-1}(F_1(b))$ is a weak homotopy equivalence for each $b\in B$.
\medskip
To apply these, let $B_n\M$ be the subspace of $B\M$ formed from the products $\Delta^p\times \M^p$ for $p\le n$, with $E_n X$ the preimage of $B_n\M$ in $EX$.  By (b) it suffices to show that $p\:E_nX\to B_n\M$ is a quasifibration, which we do by induction on $n$.  For the induction step we apply (a), decomposing $B_n\M$ into the two open sets $U$ and $V$ coming from decomposing $\Delta^n$ as the union of $\Delta^n - \bdy\Delta^n$ and an $\varepsilon$-neighborhood of $\bdy\Delta_n$ in $\Delta^n$, respectively.  Over $U$ and $U\cap V$ the space $E_nX$ is a product $X\times U$ or $X \times U\cap V$.  For $V$, a deformation retraction of the $\varepsilon$-neighborhood of $\bdy \Delta^n$ to $\bdy\Delta^n$ gives rise to a deformation $F_t$ as in (c), so by induction on $n$ the only thing to check is the condition that $F_1$ is a weak homotopy equivalence between fibers, which holds since each such map $F_1$ is $\cdot m\: X \to X$ for some $m\in\M$ and these maps are weak homotopy equivalences by hypothesis.    \qed

\proclaim Proposition D.2. If $\pi_0\M$ is a group then the natural map $\M\to\Omega B\M$ is a weak homotopy equivalence.

\pf Proof. For each $m\in\M$ the map $\cdot m\:\M\to\M$ is a homotopy equivalence with homotopy inverse $\cdot m'$ where $m'$ lies in the path-component of $\M$ giving the inverse of the path-component of $m$ with respect to the group structure on $\pi_0\M$, since homotopies from the compositions of $\cdot m$ and $\cdot m'$, in either order, to the identity are given by right multiplication by elements along paths from $m\hskip1pt m'$ and $m'm$ to $e$.  The previous lemma then implies that $\M\to E\M\to B\M$ is a cofibration.  

The contraction of $E\M$ described earlier gives a map from this cofibration to the standard path fibration $\Omega B\M \to PB\M\to B\M$, using the unique vertex of $B\M$ as the basepoint for the pathspace.  Namely, the map $E\M\to PB\M$ sends a point in $E\M$ to the path in $B\M$ which is the image in $B\M$ of the path in $E\M$ traced out by this point during the contraction of $E\M$.  In particular, a vertex $m\in E\M$ is sent to the loop in $B\M$ defined by the edge in $B\M$ labeled $m$.  This means that the restriction of the map $E\M\to PB\M$ to the fiber over the basepoint is the canonical map $\M\to\Omega B\M$.  

This map from the cofibration to the path fibration is the identity on the base spaces and is a homotopy equivalence on the total spaces since both of these are contractible.  Hence by the five lemma the map is a weak homotopy equivalence on the fibers.   \qed

The Group Completion Theorem generalizes the preceding proposition to certain cases when $\pi_0\M$ is not a group.  Initially we will assume $\pi_0\M = \Z_{\ge 0}$, so the path components $\M_p$  of $\M$ are indexed by integers $p\ge 0$.  Choose an element $m_0$ of the path component $\M_1$, so $m_0$ generates $\pi_0\M$.  Then left-multiplication by $m_0$ gives maps $m_0\cdot\:\M_p\to\M_{p+1}$, and we are interested in forming a limit object $\M_\infty$ with respect to these maps.  If the maps were suffficiently nice inclusion maps we could just take $\M_\infty = \cup_p\M_p$, but to deal with the general case we instead let $\M_\infty$ be the mapping telescope of the sequence $\M_0 \to \M_1\to \M_2\to\cdots$ where each map is $m_0\cdot$.  Recall that the mapping telescope is the union of the mapping cylinders of the maps in the sequence, and if the maps are cofibrations, the mapping telescope has the same homotopy type as $\cup_p\M_p$.  Since homology commutes with direct limits, we have $H_*(\M_\infty)=\lim_p H_*(\M_p)$.  

It will be convenient also to consider the larger mapping telescope of the sequence $\M\to\M\to\M\to\cdots$ where each map is $m_0\cdot$.  We denote this telescope by $T\M$.  Observe that $T\M\simeq\Z\times\M_\infty$ since $T\M$ is the disjoint union of a bi-infinite sequence of mapping telescopes of maps $m_0\cdot\:\M_{p}\to\M_{p+1}$.  These are the telescopes of the diagonal sequences of maps in the following diagram:

$$
\matrix{ \M_0 && \hskip-5pt\M_0 &&  \hskip-5pt\M_0  &  \cr
&\hskip-13pt\searrow\hskip-10pt && \hskip-14pt\searrow\hskip-10pt && \hskip-8pt\searrow  \cr
\M_1 &&  \hskip-5pt\M_1 &&  \hskip-5pt\M_1  & \cr
&\hskip-13pt\searrow\hskip-10pt && \hskip-14pt\searrow\hskip-10pt && \hskip-8pt\searrow \cr
\M_2 &&  \hskip-5pt\M_2 &&  \hskip-5pt\M_2  & \cr
&\hskip-13pt\searrow\hskip-10pt && \hskip-14pt\searrow\hskip-10pt && \hskip-8pt\searrow  \cr
}
$$

The monoid $\M$ acts on $E\M$ on the left, so we can form the telescope $TE\M$ of the sequence of maps $E\M\to E\M \to E\M\to \cdots$ given by left multiplication by $m_0$.  Since left and right multiplication commute, we can also form $ET\M$ using the right action of $\M$ on $T\M$, and we get the same space $ET\M=TE\M$.  Since each space $E\M$ in the telescope $TE\M$ is contractible, so is the telescope itself, and we can choose a contraction extending the earlier contraction of the initial $E\M$ in the telescope.  Having such a contraction allows us to form the following commutative diagram:

$$
\matrix{\M  &\To & E\M &\To & B\M \cr 
\downarrow &&\downarrow && \Vert \cr
T\M&\To & ET\M & \To & B\M\cr
\downarrow &&\downarrow && \Vert \cr
\Omega B\M & \To & PB\M & \To & B\M\cr
}
$$

\noindent
The top row need not be a quasifibration as it was in the previous proposition since the right-action of $\M$ on itself may not be by weak homotopy equivalences.  (The action takes $\M_p$ to $\M_{p+1}$ and these two spaces can have different fundamental groups, as happens for example for the $\M=\M^\infty$ in the proof of the Madsen-Weiss theorem.)  The middle row may not be a quasifibration either for the same sort of reason since the right-action of $\M$ on $T\M$ need not be by weak homotopy equivalences.  (The same example illustrates this.)  However, we will show that  the right-action of $\M$ on $T\M$ is by homology isomorphisms, and this will lead to the Group Completion Theorem asserting that the downward arrow in the lower left corner of the preceding diagram is a homology isomorphism.

\proclaim Group Completion Theorem. If $\pi_0\M=\Z_{\ge 0}$ and $\M$ is homotopy commutative, then the map $T\M\to \Omega B\M$ induces an isomorphism on homology.  Restricting to one component, the map $\M_\infty\to\Omega_0 B\M$ is a homology isomorphism.

\pf Proof. First we check that the right action of $\M$ on $T\M$ is by homology isomorphisms.  It suffices to show this for the right action by $m_0$.  Since $H_*(T\M)=\lim H_*(\M)$ and $\M$ is assumed to be homotopy commutative, we can instead use the left action of $m_0$ on $T\M$.  We claim that this left action is a homotopy equivalence.  This is a general fact about the action of any map $f\:X\to X$ on the mapping telescope of the sequence  $X\to X\to X \to\cdots$ where each map is $f$.  Namely, if $\sigma$ is the map of the telescope shifting each mapping cylinder to the next one, then $f\sigma =\sigma f$ is homotopic to the identity by the homotopy that slides points one unit to the right along the line segments of the mapping cylinders.

Next we have the fact that the suspension of any homology isomorphism $f\: X\to Y$ is a weak homotopy equivalence.  This can be seen as follows.  The suspension $\Sigma f\:\Sigma X\to\Sigma Y$ induces an isomorphism on $\pi_1$ since $\pi_1X$ is free with basis corresponding to all but one of the path-components $X_i$ of $X$, and similarly for $Y$, and $f$ induces a bijection $\pi_0X \to\pi_0Y$.  For the higher homotopy groups, lifting $\Sigma f$ to the universal covers gives a map $\tilde{\Sigma f}\:\tilde{\Sigma X}\to \tilde{\Sigma Y}$ that induces an isomorphism on homology since $\tilde{\Sigma X}$ consists of copies of the suspensions $\Sigma X_i$ glued together in a tree-like pattern, and similarly for $Y$ with the same tree-like pattern.  The relative Hurewicz theorem then implies that $\tilde{\Sigma f}$ induces isomorphisms on homotopy groups, hence the same is true for $\Sigma f$ itself.

Consider now the fiberwise suspension of the projection map $ET\M\to B\M$.  By definition, this is the union of two copies of the mapping cylinder of $ET\M\to B\M$ with the ends $ET\M$ identified.  Let us write the fiberwise suspension as $\Sigma_f ET\M$.  It has a projection $\Sigma_f ET\M\to B\M$ with fibers the ordinary suspension $\Sigma T\M$. Since the suspension of a homology isomorphism is a weak homotopy equivalence, Lemma~D.1 implies that we have a quasifibration $\Sigma T\M\to  \Sigma_f ET\M \to B\M$.  We also have a quasifibration $\Sigma\Omega B\M \to \Sigma_f PB\M \to B\M$ in view of the following general fact:

\proclaim Lemma D.3.  The fiberwise suspension of a fibration is a quasifibration.

Postponing the proof of this temporarily, let us finish the proof of the Group Completion Theorem.  Since fiberwise suspensions are double mapping cylinders, the map $ET\M\to PB\M$ induces a commutative diagram of quasifibrations

$$
\matrix{\Sigma T\M & \To & \Sigma_f ET\M & \To & B\M \cr
\downarrow &&\downarrow && \Vert \cr
\Sigma\Omega B\M & \To &\Sigma_f PB\M & \To & B\M \cr
}
$$

\noindent
The map between the total spaces is a weak homotopy equivalence since both total spaces are homotopy equivalent to $B\M \vee B\M$ by virtue of the contractibility of the middle spaces $ET\M$ and $PB\M$ in the two double mapping cylinders.  The five lemma applied to the long exact sequences of homotopy groups for the two quasifibrations then implies that the map of fibers $\Sigma T\M\to \Sigma\Omega B\M$ is a weak homotopy equivalence, hence an isomorphism on homology, so the map $T\M\to \Omega B\M$ is also an isomorphism on homology.  \qed

\pf Proof of Lemma D.3.  Let $p\:E\to B$ be a fibration with fiberwise suspension $\Sigma_f E\to B$.  Suppose we are given a homotopy $h_t\:X \to B$ with an initial lift $\tilde h_0\:X\to\Sigma_f E$ of $h_0$.  If $2B \subset\Sigma_f E$ is the subspace consisting of the two copies of $B$ at the target ends of the two mapping cylinders that make up $\Sigma_f E$, then $\Sigma_f E - 2B$ is a product $E\times (-1,1)$.  Part of a lifting $\tilde h_t\:X\to\Sigma_f E $ can be defined on $X - \tilde h_0^{-1}(2B)$ by taking its first coordinate in $E$ to be a lift of $h_t$ to $E$ and letting its second coordinate in $(-1,1)$ be independent of $t$.  On the rest of $X$ we could let $\tilde h_t$ be just $h_t$ mapping to the appropriate copy of $B$ in $2B$.  However, it is not clear that the resulting $\tilde h_t$ is continuous on all of $X$ since we have defined it on two subspaces, one of which is open and the other is closed, rather than on two subspaces which are both open or both closed.

To get around this problem we use the elementary fact that in order to have a quasifibration it is not necessary to lift the homotopy $h_t$  exactly,  it suffices just to lift a modification of $h_t$ obtained by reparametrizing the $t$ interval, namely by first letting $h_t$ be stationary for $t$ in $[0,1/2]$, then letting it do what it did before but twice as fast so that it ends with the same $h_1$.  To construct the lifting $\tilde h_t$ of the new $h_t$, we first let $\tilde h_t=r_t \tilde h_0$ during the $t$ interval $[0,1/2]$, where $r_t\:\Sigma_f E\to\Sigma_f E$ is a homotopy that is fixed on $2B$ and gradually expands the interval $[-1/2,1/2]$ to $[-1,1]$ in the second factor of $E\times (-1,1)$.  Then for $t$ in $[1/2,1]$ we take the composition $r_1\tilde h_t$ using the $\tilde h_t$ constructed in the previous paragraph, with the $t$ interval reparametrized from $[0,1]$ to $[1/2,1]$.  With these modifications the new $\tilde h_t$ is continuous since it is continuous on the closed sets $\tilde h_0^{-1}(E\times [-1/2,1/2])$ and $X-\tilde h_0^{-1}(E\times (-1/2,1/2))$.  \qed

\subsect A Modest Generalization

The condition that $\pi_0\M$ is $\Z_{\ge 0}$ can be weakened to the assumption that $\pi_0\M$ is just finitely generated.  In this case, let $m_1,\cdots,m_k$ be a set of generators and let $m_0$ be the product $m_1 m_2 \cdots m_k$.  With this $m_0$ we can again form the telescope $T\M$.  As before, the right action of $m_0$ on $T\M$ is by homology isomorphisms.   This implies the right action by any of the generators $m_i$ is by homology isomorphisms since the fact that $m_1m_2\cdots m_k$ acts as homology isomorphisms implies that $m_1$ acts injectively on homology and $m_k$ acts surjectively, but by homotopy commutativity this implies all the generators act by homology isomorphisms.  It follows that the right action by all elements of $\M$ on $T\M$ is by homology isomorphisms.

The rest of the proof goes through unchanged to give an isomorphism $H_*(T\M)\approx H_*( \Omega B\M)$.  One can check that all the path-components of $T\M$ are homotopy equivalent.  
Denoting this homotopy type by $\M_\infty$ (it does not depend on the choice of generators), we then obtain an isomorphism $H_*(\M_\infty) \approx H_*(\Omega_0 B\M)$.  

With an appropriate reformulation of the theorem, the finite generation assumption on $\pi_0\M$ can in fact be omitted.  See [MS].

\leftskip 20pt
\parindent=-20pt

\vskip-15pt

\section References

[A] J.F. Adams, Infinite Loop Spaces, Annals of Math. Studies 90, Princeton Univ. Press, 1978.

[BP] M. Barratt and S. Priddy, On the homology of non-connected monoids and their associated groups, Comment. Math. Helv. 47 (1972), 1--14.

[BT] C.F. B\"odigheimer and U. Tillmann, Stripping and splitting decorated mapping class groups,  in Cohomological methods in homotopy theory (Bellaterra, 1998), Progr. Math. 196 (2001),
47--57.

[B] S.K. Boldsen, Improved homological stability for the mapping class group with integral or twisted coefficients,  Math. Z. 270 (2012), 297--329.

[C1] J. Cerf, Th\'eor\`emes de fibration des espaces de plongements. Applications.  S\'em. H. Cartan, v. 15 no. 8 (1962-63), 1-13.

[C2] J. Cerf, Sur les diff\'eomorphismes de la sph\`ere de dimension trois ($\Gamma_4=0$), Springer Lecture Notes 53 (1968).

[EE] C. J. Earle and J. Eells, A fibre bundle description of Teichm\" uller theory, J. Differential Geometry 3 (1969), 19--43.

[ES] C. J. Earle and A. Schatz, Teichm\" uller theory for surfaces with boundary,  
J. Differential Geometry 4 (1970), 169--185.

[G1] S. Galatius, Stable homology of automorphism groups of free groups, Ann. of Math. 173 (2011), 705--768.

[G2] S. Galatius, Mod p homology of the stable mapping class group, Topology 43 (2004), 1105--1132.

[GMTW] S. Galatius, I. Madsen, U. Tillmann, and M. Weiss, The homotopy type of the cobordism category, Acta Math. 202 (2009), 195--239.

[GRW] S. Galatius and O. Randal-Williams, Monoids of moduli spaces of manifolds, Geom. Topol. 14 (2010), 1243--1302.

[Gr] A. Gramain, Le type d'homotopie du groupe des diff\'eomorphismes d'une surface compacte, Ann. Scient. \'Ec. Norm. Sup. 6 (1973), 53--66.

[Han] E. Hanbury, Homological stability of non-orientable mapping class groups with marked points,  Proc. A.M.S. 137 (2009),  385--392.

[Har] J.L. Harer, Stability of the homology of the mapping class groups of orientable surfaces, Ann. of Math. 121 (1985), 215--249.

[H1] A. Hatcher, {\it Algebraic Topology}, Cambridge Univ. Press 2002.

[H2] A. Hatcher, {\it Vector Bundles and K-Theory}, notes available from the author's webpage.

[HW] A. Hatcher and N. Wahl, Stabilization for mapping class groups of 3-manifolds, Duke Math. J. 155 (2010), 205--269. arXiv:0709.2173

[I] N.V. Ivanov, Stabilization of the homology of Teichm\"uller modular groups, Algebra i Analiz 1 (1989), 110--126 (Russian); translation in Leningrad Math. J. 1 (1990), 675--691.

[MT] I. Madsen and U. Tillmann, The stable mapping class group and $Q(\hbox{\Bbb C}P^\infty_+)$, Invent. Math. 145 (2001), 509--544.

[MW] I. Madsen and M. Weiss, The stable moduli space of Riemann surfaces: Mumford's conjecture, Ann. of Math. 165 (2007), 843--941. 

[M] D. McDuff, Configuration spaces of positive and negative particles, Topology 14 (1975), 91--107.

[MS] D. McDuff and G. Segal, Homology fibrations and the ``group-completion" theorem, Invent. Math. 31 (1976), 279--284.

[RW] O. Randal-Williams, Resolutions of moduli spaces and homological stability, \hfill\break arXiv:0909.4278. 

[Sa] N.N. Savelev, Quillenization of classifying spaces of infinite symmetric powers of discrete groups, Vestnik Moskov. Univ. Ser. I Mat. Mekh. (1989), no.3, 72--74 (Russian); translation in Moscow Univ. Math. Bull. 44 (1989), no. 3, 74--76.

[S1] G. Segal, Configuration-spaces and iterated loop-spaces, Invent. Math. 21 (1973), 213--221. 

[S2] G. Segal, The topology of spaces of rational functions, Acta Math. 143 (1979), 39--72.

[T] U. Tillmann, On the homotopy of the stable mapping class group. Invent. Math. 130 (1997), 257--275.

[W1] N. Wahl, Homological stability for the mapping class groups of non-orientable surfaces, Invent. Math. 171 (2008), 389--424.

[W2] N. Wahl, Homological stability for mapping class groups of surfaces, Handbook of Moduli, Vol. III, Advanced Lectures in Mathematics 26 (2012), 547--583.

\end